\documentclass[12pt,leqno,a4paper]{article}

\usepackage{amsfonts}
\usepackage{latexsym}
\usepackage{amssymb,amsmath}
\usepackage{epsfig}
\usepackage{hyperref}
\usepackage{color}
\everymath={\displaystyle}

\newcommand{\R}{\mathbb R}

\newcommand{\be}{\begin{equation}}
\newcommand{\ee}{\end{equation}}
\newcommand{\ba}{\begin{eqnarray}}
\newcommand{\ea}{\end{eqnarray}}

\newcommand{\tr}{^\text{tr}}

\newcommand{\cqfd}
{%
\mbox{}%
\nolinebreak%
\hfill%
\rule{2mm}{2mm}%
\medbreak%
\par%
}

\numberwithin{equation}{section}

\newtheorem{thm}{\bf Theorem}
\newtheorem{lem}[thm]{\bf Lemma}

\newtheorem{prop}[thm]{\bf Proposition}

\newtheorem{rem}[thm]{\bf Remark}

\renewcommand{\leq}{\leqslant}

\date{\empty}

\title{Control of three heat equations coupled with two cubic nonlinearities}
\author{Jean-Michel
Coron\thanks{Universit\'{e} Pierre et
Marie Curie-Paris 6, UMR 7598 Laboratoire Jacques-Louis Lions, 75005 Paris,
France. E-mail: \texttt{coron@ann.jussieu.fr}. JMC was supported by the ERC
advanced grant 266907 (CPDENL) of the 7th Research Framework Programme
(FP7).}, Jean-Philippe Guilleron\thanks{Universit\'{e} Pierre et
Marie Curie-Paris 6, UMR 7598 Laboratoire Jacques-Louis Lions, 75005 Paris,
France. E-mail: guilleron@ann.jussieu.fr.
JPG was partially supported by the ERC
advanced grant 266907 (CPDENL) of the 7th Research Framework Programme
(FP7).}}

\begin{document}
\maketitle
\begin{abstract}
We study the null controllability of three parabolic equations. The control is acting only on one of the three
equations. The three equations are coupled by means of two cubic nonlinearities. The linearized control system around $0$ is not null controllable. However, using the cubic nonlinearities, we prove the (global) null controllability of the control system. The proof relies on the return method, an algebraic solvability and smoothing properties of the parabolic equations.
\end{abstract}
\section{Introduction}
\setcounter{equation}{0}
Let $N$ be a positive integer and let $\Omega$ be a nonempty connected bounded subset of $\R^N$ of class $C^2$. Let
 $\omega$ be a nonempty open subset of $\Omega$. We denote by $\chi_\omega:\Omega \to \R$  the characteristic function of $\omega$ and let $T\in (0,+\infty)$. We are interested in the control system
\begin{equation}
\label{eqsys}
\left\{
\begin{array}{ll}
\alpha_t-\Delta \alpha =\beta^3 & \text{ in }(0,T)\times  \Omega,
\\
\beta_t-\Delta \beta =\gamma^3 & \text{ in }(0,T)\times  \Omega,
\\
\gamma_t-\Delta \gamma= u\chi_\omega& \text{ in }(0,T)\times  \Omega,
\\
\alpha=\beta=\gamma=0 &\text{ in }(0,T)\times \partial \Omega.
\end{array}
\right.
\end{equation}
It is a control system where, at time $t\in [0,T]$,
the state is $(\alpha(t,\cdot),\beta(t,\cdot),\gamma(t,\cdot))\tr:
\Omega\to \R^3$ and the control is $u(t,\cdot):\Omega\to \R$. Let us point out that, due to the recursive structure of \eqref{eqsys} (one first solves the last parabolic equation of \eqref{eqsys}, then the second one and finally the first one), it follows from classical results on linear parabolic equations that the  Cauchy problem associated to \eqref{eqsys} is globally  well-posed in the $L^\infty$ setting, i.e. with bounded measurable initial data, controls, and solutions.

The main goal of this paper is to prove the following global null controllability
result on control system~\eqref{eqsys}.
\begin{thm}
\label{thnullcont}
 For every $(\alpha^0,\beta^0,\gamma^0)\tr \in L^\infty(\Omega)^3$,
there exists a control $u\in L^\infty((0,T)\times\Omega)$ such that the
solution $(\alpha,\beta,\gamma)\tr\in L^\infty((0,T)\times\Omega)^3$ to the Cauchy problem
\begin{equation}
\label{Cauchy}
\left\{
\begin{array}{ll}
\alpha_t-\Delta \alpha =\beta^3 & \text{ in }(0,T)\times  \Omega,
\\
\beta_t-\Delta \beta =\gamma^3 & \text{ in }(0,T)\times  \Omega,
\\
\gamma_t-\Delta \gamma= u\chi_\omega& \text{ in }(0,T)\times  \Omega,
\\
\alpha=\beta=\gamma=0 &\text{ in }(0,T)\times \partial \Omega,
\\
\alpha(0,\cdot)=\alpha^0(\cdot), \, \beta(0,\cdot)=\beta^0(\cdot), \,
\gamma(0,\cdot)=\gamma^0(\cdot) & \text{ in } \Omega,
\end{array}
\right.
\end{equation}
satisfies
\begin{equation}\label{finalnul}
  \alpha(T,\cdot)=\beta(T,\cdot)=\gamma(T,\cdot)=0 \text{ in }\Omega.
\end{equation}
\end{thm}

The controllability of systems of partial differential equations with a small number of controls is an important subject which has been
recently investigated in a large number of articles. For the case of linear systems, let us mention in particular
\begin{itemize}
\item For systems of parabolic equations in dimension 1 or larger:
\cite{2009-de-Teresa-Zuazua, 2000-deTeresa-CPDE, 2006-Gonzalez-Burgos-Perez-Garcia-AA, 2007-Guerrero-SICON}.  A key step in these papers is to establish suitable Carleman estimates. In dimension 1, the method of moments can lead to very precise (and sometimes unexpected) results; see, in particular
\cite{2014-Ammar-Khodja-Benabdallah-Gonzalez-Burgos-de-Teresa-CRAS, 2014-Ammar-Khodja-Benabdallah-Gonzalez-Burgos-de-Teresa-JFA,
2014-Benabdallah-Boyer-Gonzalez-Burgos-Olive-SICON, 2014-Boyer-Olive-MCRF}. See also the survey paper \cite{2011-Ammar-Khodja-Benabdallah-Gonzalez-Burgos-de-Teresa-MCRF} and the reference therein.
\item For systems of Schr\"{o}dinger equations: \cite{2014-Alabau-MCSS}, which uses transmutation together with a controllability result for systems of
wave equations proved in the same article. See also \cite{2011-Rosier-de-Teresa-CRAS} for the controllability of a cascade system of conservative equations.
\item For Stokes equations of incompressible fluids:
\cite{2006-Fernandez-Cara-Guerrero-Imanuvilov-Puel-SICON, 2007-Guerrero-AIHP, 2009-Coron-Guerrero-JMPA, 2013-Carreno-Guerrero-JMFM}.
Again Carleman estimates are key ingredients here.
\item For hyperbolic equations: \cite{2013-Alabau-ADE, 2014-Alabau-MCSS}, which rely on multiplier methods, and \cite{2013-Alabau-Leautaud-JMPA} which uses microlocal analysis.
\end{itemize}
Let us assume that $0$ is a trajectory (i.e. a solution) of the system of partial differential equations.
If the linearized control system is controllable, one can expect to get the local null controllability.
For systems of partial differential equations with a small number of controls it has been proven to be the case,
for example, for the Navier Stokes equations in \cite{2013-Carreno-Guerrero-JMFM}.

Note that the linearized control system of \eqref{eqsys} around $0$ is clearly not controllable. When  the linearized control system around $0$ is not controllable one may still expect that
the nonlinearities can give the controllability. A method to treat this case is the return method. It consists in looking for (nonzero) trajectories of the control system going from $0$ to $0$ such that the linearized control system is controllable. This method has been introduced in \cite{1992-Coron-MCSS} for a stabilization issue and used for the first time in \cite{1996-Coron-JMPA} to
get the controllability of a partial differential equation (the Euler equation of incompressible fluids). This method can also
be used to get controllability of systems of partial differential equations with a small number of controls. See, for example,
\begin{itemize}
\item \cite{2002-Coron-COCV} for a water tank control system modeled by means of the Saint-Venant equations.
\item \cite{2009-Coron-Guerrero-JMPA, 2014-Coron-Lissy-IM} for the Navier-Stokes equations.
\item \cite{2010-Coron-Guerrero-Rosier-SICON} for a system of two nonlinear heat equations.
\end{itemize}
Let us give more details about \cite{2010-Coron-Guerrero-Rosier-SICON} since it deals with a control system related to our system \eqref{eqsys}. The control system considered in \cite{2010-Coron-Guerrero-Rosier-SICON} is
\begin{equation}
\label{eqsys-old}
\left\{
\begin{array}{ll}
\beta_t-\Delta \beta =\gamma^3 & \text{ in }(0,T)\times  \Omega,
\\
\gamma_t-\Delta \gamma= u\chi_\omega& \text{ in }(0,T)\times  \Omega,
\\
\beta=\gamma=0 &\text{ in }(0,T)\times \partial \Omega,
\end{array}
\right.
\end{equation}
where, at time $t\in [0,T]$,
the state is $(\beta(t,\cdot),\gamma(t,\cdot))\tr:
\Omega\to \R^2$ and the control is $u(t,\cdot):\Omega\to \R$. (In fact, slightly more general control systems of two coupled parabolic equations are considered in \cite{2010-Coron-Guerrero-Rosier-SICON}.) Using the return method, it is proved in
\cite{2010-Coron-Guerrero-Rosier-SICON} that the control system \eqref{eqsys-old} is locally null
controllable. We use the same method here. However the construction of trajectories of the control system going from $0$ to $0$ such that the linearized control system is (null) controllable is much more complicated for
the control system \eqref{eqsys} than for the control system \eqref{eqsys-old}.

The construction of trajectories of the control system \eqref{eqsys} going
from $0$ to $0$ such that the linearized control system is (null)
controllable follows from simple scaling arguments (see \eqref{defbaralpha} to  \eqref{defbaru} below) and the  following theorem.
\begin{thm}\label{thmexistenceabc}
There exists $(a,b,c)\in C_0^{\infty}(\mathbb{R}\times\mathbb{R})^3$ such that
\begin{gather}
\label{supportabc}
\text{the supports of $a$, $b$, and $c$ are included in $[-1,1]\times [-1,1]$},
\\
\label{nonzeroensemble}
\{(t,r); \, r>0, \, b(t,r)\not= 0 \textrm{ and } c(t,r)\not= 0\}\not =\emptyset,
\\
\label{abceven}
a(t,r)=a(t,-r), \, b(t,r)=b(t,-r), \, c(t,r)=c(t,-r), \, \forall (t,r)\in \R\times\R,
\\
\label{edpab3}
a_t-a_{rr}-\frac{N-1}{r}a_{r} =b^3 \text{ in }  \mathbb{R}\times\mathbb{R}^*,
\\
\label{edpbc3}
b_t-b_{rr}-\frac{N-1}{r}b_{r} =c^3 \text{ in } \mathbb{R}\times\mathbb{R}^*.
\end{gather}
\end{thm}

An important ingredient of the proof of Theorem~\ref{thmexistenceabc} is the
following proposition which is related to Theorem~\ref{thmexistenceabc} in the
stationary case.
\begin{prop}
\label{thmexistenceabc-stationary}
 There exists $(A,B,C)\in C^{\infty}(\R)^3$ and $\delta_A\in (0,1/2)$  such that
\begin{gather}
\label{supportABCgood}
\text{ the supports of $A$, $B$, and $C$ are included in $[-1,1]$},
\\
\label{BCnotempty}
\{z; \,  z>0, B(z)\not= 0 \text{ and } C(z)\not= 0\}\not =\emptyset,
\\
\label{ABCeven}
A(z)=A(-z), \, B(z)=B(-z), \, C(z)=C(-z), \, \forall z\in \R,
\\
\label{barApourprochede1<1-th}
{ A}(z)= e^{-1/(1-z^2)}
\text{ if }\ 1-\delta_A < z <1,
\\
\label{eqAB3}
-A''-\frac{N-1}{z}A'=B^3 \text{ in } \mathbb{R}^*,
\\
\label{eqBC3}
-B''-\frac{N-1}{z}B'=C^3 \text{ in } \mathbb{R}^*,
\\
\label{eqBnulseulementen1/2}
\left(B(z)=0 \text{ and } z\in[0,1)\right)\Leftrightarrow \left(z=\frac{1}{2}\right),
\\
\label{deriveeB1/2<0}
B'\left(\frac{1}{2}\right)<0,
\\
\label{C(1/2)notzero}
C\left(\frac{1}{2}\right) > 0,
\\
\label{C'enzerodeC}
\left(C(z)=0 \text{ and } z \in [0,1)\right)\Rightarrow
\left( z \in (0,1)\text{ and } C'(z)\not =0\right).
\end{gather}
\end{prop}
This proposition is proved in Section~\ref{secstationary}. In Section~\ref{sectimevarying}
we show how to use Proposition~\ref{thmexistenceabc-stationary} in order to
prove Theorem~\ref{thmexistenceabc}. Finally, in Section~\ref{sec-proof-local-controllability}, we deduce
Theorem~\ref{thnullcont} from
Theorem~\ref{thmexistenceabc}.

\begin{rem}\label{rem-extension}
Looking to our proof of Theorem~\ref{thnullcont}, it is natural to conjecture that this theorem still holds if, in \eqref{Cauchy}, $\beta^3$ and $\gamma^3$ are replaced
by $\beta^{2p+1}$ and $\gamma^{2q+1}$ respectively, where $p$ and $q$ are arbitrary nonnegative integers.
\end{rem}

\section{Proof of Proposition~\ref{thmexistenceabc-stationary} (stationary case)}
\label{secstationary}
\setcounter{equation}{0}
In order to construct $A$, one shall use the following lemma.

\begin{lem}
\label{lemmabonG}
There exists $\delta_0\in (0,1)$ such that, for every $\delta \in (0,\delta_0)$,
there exists a function $G\in C^\infty ([0,+\infty ))$
such that
\begin{gather}
\label{Gpourzprochede1/2}
G(z) = z^3 \left(z-\frac{1}{2}\right)^3 \text{ for }
\frac{1}{2} -\delta <z<\frac{1}{2} +\delta,
\\
\label{signeG}
(z-\frac{1}{2})G(z)>0  \text{ for } 0<z<1,\ z\ne \frac{1}{2},
\\
\label{finiteptpbforG}
\left\{z\in (0,1);\, (G^{1/3})''(z)+\frac{N-1}{z}(G^{1/3})'(z)=0\right\} \text{ is finite},
\end{gather}
and such that the solution ${ A}:(0,+\infty)\rightarrow \R$ to the Cauchy problem
\begin{gather}
\label{eqCauchyz}
{ A}(1)={ A}'(1)=0,\, { A}''(z)+\frac{N-1}{z}{ A}'(z) = G(z),\, z>0,
\end{gather}
satisfies
\begin{gather}
\label{barApourprochede0}
\text{there exists $c_0\in \R$ such that } { A}(z)= c_0-z^{8}\text{ if }\  0<z<\delta,
\\
\label{barApourprochede1<1}
{ A}(z)= e^{-1/(1-z^2)}
\text{ if }\ 1-\delta < z <1,
\\
\label{barApourz>1}
{ A}(z)= 0\text{ if } z\in  [1,+\infty).
\end{gather}
\end{lem}

\textbf{Proof of Lemma~\ref{lemmabonG}.} Let us first emphasize that it follows from \eqref{Gpourzprochede1/2} and \eqref{signeG} that
$G^{1/3}$ is of class $C^\infty$ on $(0,1)$, hence \eqref{finiteptpbforG} makes sense.
Let $\delta \in(0,1/4)$. Let ${\bar G}\in C^\infty([0,+\infty))$ be such that \eqref{Gpourzprochede1/2} and \eqref{signeG} hold for $G={\bar G}$ and
\begin{gather}
\label{valueG0petit}
{\bar G}(z)=-8(6+N)z^6, \forall z\in (0,\delta),
\\
\label{valueG0zproche1}
{\bar G}(z)= \left(\frac{-2+6z^4}{(1-z^2)^4}-\frac{2(N-1)}{(1-z^2)^2}\right)e^{-1/(1-z^2)},\,
\forall z\in ((1-\delta), 1),
\\
\label{valueG0z>1}
{\bar G}(z)=0,\,  \forall z\in (1,+\infty),
\\
\label{G0analytic}
{\bar G} \text{ is analytic on } (0,1)\setminus \{\delta, (1/2)-\delta, (1/2)+\delta, 1-\delta\}.
\end{gather}
One easily sees that such $\bar G$ exists if $\delta \in (0,1/4)$ is small enough, the smallness depending on $N$.
Frow now on, $\delta$ is always assumed to be small enough. Let $\kappa\in \R$. Let us define $G\in C^\infty([0,+\infty))$  by
\begin{gather}
\label{Galpha=G0}
G:={\bar G} \text{ in }[0,\delta]\cup [(1/2)-\delta, (1/2)+\delta]\cup[1-\delta,+\infty),
\\
\label{valueGdeuxiemeinterval}
G(z):={\bar G}(z)+\min\{\kappa,0\}e^{-1/(z-\delta)}e^{-1/(1-2\delta-2z)},\, \forall z \in (\delta, (1/2)-\delta),
\\
\label{valueGquatriemeinterval}
G(z):={\bar G}(z)+\max\{\kappa,0\}e^{-1/(2z-1-2\delta)}e^{-1/(1-\delta-z)},\, \forall z \in ((1/2)+\delta,1-\delta).
\end{gather}
Let ${ A}$ be the solution of the Cauchy problem \eqref{eqCauchyz}. From \eqref{Galpha=G0}, one has \eqref{Gpourzprochede1/2} and \eqref{signeG}. From \eqref{G0analytic}, \eqref{valueGdeuxiemeinterval}, and \eqref{valueGquatriemeinterval}, one gets that
\begin{equation}\label{Ganalytic}
G \text{ is analytic on } (0,1)\setminus \{\delta, (1/2)-\delta, (1/2)+\delta, 1-\delta\},
\end{equation}
which implies \eqref{finiteptpbforG} since $(G^{1/3})''$ cannot be identically equal to $0$ on one of the
 five intervals $(0,\delta)$, $(\delta,(1/2)-\delta)$, $((1/2)-\delta,(1/2)+\delta)$, $((1/2)+\delta,1-\delta)$, and $(1-\delta,1)$.
\begin{rem}
We require \eqref{Ganalytic} only to get \eqref{finiteptpbforG}. However \eqref{finiteptpbforG}
 can  also be obtained without requiring \eqref{Ganalytic} by using
genericity arguments.
\end{rem}
 From \eqref{eqCauchyz}, \eqref{valueG0zproche1}, and \eqref{Galpha=G0}, one gets \eqref{barApourprochede1<1}. From \eqref{eqCauchyz}, \eqref{valueG0z>1}, and \eqref{Galpha=G0}, one gets \eqref{barApourz>1}.

It remains to prove that, for some $\kappa\in \R$, one has \eqref{barApourprochede0}. Let us first point out
that, for every $y\in C^2((0,\delta))$,
\begin{multline}\label{solpart}
  \left(y''+\frac{N-1}{z}y'=0\right)\Rightarrow
  \\
  \left(\exists \, (c_0,c_1)\in \R^2
  \text{ such that } y(z)=c_0+c_1E(z),\, \forall z\in (0,\delta)\right),
\end{multline}
where
\begin{gather}
\label{defEN>2}
\text{if } N\not =2,\,  E(z):=\frac{1}{(2-N)z^{N-2}}, \, \forall z\in(0,+\infty),
\\
\label{defEN=2}
\text{if } N=2,\,  E(z):=-\ln(z), \, \forall z\in(0,+\infty).
\end{gather}
From \eqref{eqCauchyz}, \eqref{valueG0petit}, \eqref{Galpha=G0}, one gets that $y:={ A}+z^8$ satisfies the assumption
 of  the implication \eqref{solpart}. Hence, by  \eqref{solpart}, one gets the existence of $(c_0,c_1)\in \R^2$ such that
\begin{equation}\label{valuebarASC1}
  { A}(z)= c_0-z^8+c_1E(z), \, \forall z\in (0,\delta).
\end{equation}
It suffices to check that, for some $\kappa \in \R$,
\begin{equation}\label{C1=0}
  c_1=0.
\end{equation}
From \eqref{eqCauchyz}, one has
\begin{gather}
\label{expressionbarAN>2}
  \text{if $N\not = 2$, }{ A}(z)=-\frac{1}{(N-2)z^{N-2}}\int_1^z s^{N-1}G(s)ds+\frac{1}{N-2} \int_1^z s G(s)ds, \, \forall z \in (0,1],
\\
\label{expressionbarAN22}
\text{if $N=2$, }{ A}(z)=\ln(z)\int_1^z sG(s)ds- \int_1^z s \ln(s) G(s)ds, \, \forall z \in (0,1],
\end{gather}
which, together with \eqref{defEN>2}, \eqref{defEN=2}, \eqref{valuebarASC1}, with $z\rightarrow0$, gives
\begin{equation}\label{valueC1}
  c_1=\int_0^1s^{N-1}G(s)ds.
\end{equation}
 From \eqref{Galpha=G0}, \eqref{valueGdeuxiemeinterval}, and \eqref{valueGquatriemeinterval}, one has
\begin{equation}\label{limitalpha}
  \lim_{\kappa \rightarrow +\infty} \int_0^1s^{N-1}G(s)ds=+\infty \text{ and }
  \lim_{\kappa \rightarrow -\infty} \int_0^1s^{N-1}G(s)ds=-\infty.
\end{equation}
In particular, with the intermediate value theorem, there exists $\kappa \in \R$ such that
\begin{equation}\label{int=0}
  \int_0^1s^{N-1}G(s)ds=0,
\end{equation}
which, together with \eqref{valueC1},  concludes the proof of Lemma~\ref{lemmabonG}.
\cqfd

We go back to the proof of Proposition~\ref{thmexistenceabc-stationary}. We extend $ A$ to all of $\R$ by
requiring
\begin{gather}\label{extbarAen0}
   A(0)=c_0, \\
  \label{extbarAznegatif}
   A(z)=A(-z), \, \forall z\in (-\infty,0).
\end{gather}
By \eqref{barApourprochede0}, \eqref{extbarAen0}, and \eqref{extbarAznegatif}, $ A \in C^\infty(\R)$. Let
$ B\in C^0(\R^*)$ be defined by
\begin{equation}\label{defbarB}
   B:=-\left( A''+\frac{N-1}{z} A'\right)^{1/3}.
\end{equation}
 From \eqref{extbarAznegatif} and \eqref{defbarB}, one gets that
\begin{equation}\label{Bpair}
   B(z)= B(-z),\, \forall z \in \R^*.
\end{equation}
 From \eqref{defbarB}, one sees that
\begin{equation}\label{lieubarBCinfty}
   B \text{ is of class $C^\infty$ on the set $\{z\in \R^*; \, B(z)\not =0 \}$}.
\end{equation}
 From \eqref{barApourprochede0}, \eqref{extbarAznegatif}, and \eqref{defbarB}, one has
\begin{equation}\label{valuebarBprochede0}
   B(z)=2(6+N)^{1/3}z^2, \forall z\in (-\delta,\delta)\setminus\{0\},
\end{equation}
which allows to extend $ B$ to all of $\R$ by continuity by requiring
\begin{equation}\label{valuebarBen0}
   B(0)=0.
\end{equation}
 From \eqref{valuebarBprochede0} and \eqref{valuebarBen0}, we get that
\begin{equation}\label{barBCinftypresde0}
   B \text{ is of class $C^\infty$ in $(-\delta,\delta)$}.
\end{equation}
 From \eqref{signeG}, \eqref{eqCauchyz}, and \eqref{defbarB}, one gets
that
\begin{equation}\label{barBpositif}
   B \not = 0 \text{ in } (0,1)\setminus \{1/2\},
\end{equation}
which, with \eqref{lieubarBCinfty}, implies that
\begin{equation}\label{barBCinftysur01-1/2}
   B \text{ is of class $C^\infty$ in $(0,1)\setminus \{1/2\}$}.
\end{equation}
 From \eqref{Gpourzprochede1/2}, \eqref{eqCauchyz}, and \eqref{defbarB}, one has
\begin{equation}\label{valuebarBproche1/2}
   B(z)=-z\left(z-\frac{1}{2}\right), \, \forall z \in \left(\frac{1}{2}-\delta, \frac{1}{2}+\delta\right).
\end{equation}
In particular \eqref{deriveeB1/2<0} holds.  From \eqref{barApourprochede1<1} and \eqref{defbarB}, one gets
\begin{equation}\label{valeurBprochede1<1}
   B(z)=-\left(\frac{-2+6z^4}{(1-z^2)^4}-\frac{2(N-1)}{(1-z^2)^2}\right)^{1/3}e^{-1/(3-3z^2)},\,
\forall z \in (1-\delta,1),
\end{equation}
which implies the existence of $\delta_0>0$  such that, for every $\delta \in (0,
\delta_0]$,
\begin{equation}\label{barBpositifpresde1<1}
   B<0 \text{ in } (1-\delta, 1).
\end{equation}
 From \eqref{barApourz>1} and \eqref{defbarB}, one gets
\begin{equation}\label{barBpourz>1}
   B(z)=0, \forall z\in (1,+\infty),
\end{equation}
which, together with \eqref{valeurBprochede1<1}, implies that
\begin{equation}\label{barBCinfty>1}
   B \text{ is of class $C^\infty$ in $(1-\delta,+\infty)$}.
\end{equation}
 From  \eqref{Bpair}, \eqref{barBCinftypresde0}, \eqref{barBCinftysur01-1/2}, \eqref{valuebarBproche1/2},
 and \eqref{barBCinfty>1}, one gets that
\begin{equation}\label{barBpartout}
   B \text{ is of class $C^\infty$ in $\R$.}
\end{equation}

Let us now define $ C\in C^0(\R^*)$ by
\begin{equation}\label{defbarC}
   C(z):= -\left(B''+\frac{N-1}{z}B'\right)^{1/3},\, \forall z \in \R^*.
\end{equation}
 From \eqref{Bpair} and \eqref{defbarC}, one has
\begin{equation}\label{Cpair}
   C(z)= C(-z),\, \forall z \in \R^*.
\end{equation}
 From \eqref{barBpartout} and \eqref{defbarC}, one gets that
\begin{equation}\label{lieubarCCinfty}
   C \text{ is of class $C^\infty$ on the set $\{z\in \R^*; \, C(z)\not =0 \}$}.
\end{equation}
 From \eqref{valuebarBprochede0} and \eqref{defbarC}, one has
\begin{equation}\label{valuebarCprochede0}
  C(z)=-(4N)^\frac{1}{3}(6+N)^\frac{1}{9}<0, \forall z\in [-\delta,\delta].
\end{equation}
 From \eqref{valuebarBproche1/2} and \eqref{defbarC}, one has
\begin{equation}\label{valuebarCproche1/2}
   C(z)=\left(2N-\frac{N-1}{2z}\right)^{1/3}, \, \forall z \in \left[\frac{1}{2}-\delta, \frac{1}{2}+\delta\right].
\end{equation}
In particular, since $\delta>0$ is small enough,
\begin{equation}\label{Cpresde1sur2}
  C \text{ is positive and of class $C^\infty$ on $\left[\frac{1}{2}-\delta, \frac{1}{2}+\delta\right]$}.
\end{equation}
 From \eqref{valeurBprochede1<1}, \eqref{barBpourz>1}, and \eqref{defbarC}, one gets that
\begin{equation}\label{barcpresde1}
   C>0 \text{ in } [1-\delta, 1) \text{ and $ C$ is of class $C^\infty$ in $[1-\delta,+\infty)$.}
\end{equation}
 From \eqref{Cpair}, \eqref{lieubarCCinfty}, \eqref {valuebarCprochede0},  \eqref{Cpresde1sur2}, and
 \eqref{barcpresde1}, one sees that
\begin{equation}\label{barCofclassCinfty}
   C \in C^\infty(\R)
\end{equation}
if
\begin{equation}\label{goodzeroes}
C \text{ is of class } C^\infty \text{ in } (\delta,(1/2)-\delta)\cup ((1/2)+\delta,1-\delta).
\end{equation}
Let us first point out that,
by \eqref{finiteptpbforG}, \eqref{eqCauchyz}, \eqref{defbarB}, and \eqref{defbarC},
\begin{equation}\label{mauvaisfinite}
  \text{the set of  $z_0\in (\delta,(1/2)-\delta)\cup ((1/2)+\delta,1-\delta)$
such that $C(z_0)=0$ is finite.}
\end{equation}
We are going to prove that \eqref{goodzeroes} indeed holds provided that one no longer requires \eqref{Ganalytic} and that  one modifies $G$
 in a neighborhood of every
 $z_0\in (\delta,(1/2)-\delta)\cup ((1/2)+\delta,1-\delta)$ such that $C(z_0)=0$. Since $G=-B^3$,
 this comes from the following lemma.

\begin{lem}
\label{modifBG}
Let $\nu>0$, $\zeta>0$, and $\eta >0$ be such that
$[\zeta-\eta,\zeta+\eta]\subset (0,+\infty)$.
Let $B\in C^\infty ([\zeta-\eta,\zeta+\eta])$ be such that
\begin{gather}\label{propSignG''}
  B''(z)+\frac{N-1}{z}B'(z)\not =0, \, \forall z \in [\zeta-\eta,\zeta+\eta]\setminus\{\zeta\}.
\end{gather}
Then, there exists $\bar B \in C^\infty ([\zeta-\eta,\zeta+\eta])$ satisfying
\begin{gather}\label{barBprocheB}
  |\bar B (z)-B(z)|\leqslant \nu,\, \forall z\in [\zeta-\eta,\zeta+\eta],  \\
  \label{supportbarB}
  \text{the support of $\bar B-B$ is included in $(\zeta-\eta,\zeta+\eta)$,}\\
  \label{barB''1/3}
  \left(\bar B''+\frac{N-1}{z}\bar B'\right)^{1/3} \in C^\infty([\zeta-\eta,\zeta+\eta])
\end{gather}
and such that, if $\bar A \in C^\infty ([\zeta-\eta,\zeta+\eta])$ is the solution of
\begin{gather}\label{equationbarA}
  \bar A''+\frac{N-1}{z}\bar A' =-\bar B^3, \\
  \label{conditioninitbarA}
  \bar A(\zeta-\eta)=A(\zeta-\eta),\,\bar A'(\zeta-\eta)=A'(\zeta-\eta),
\end{gather}
then,
\begin{equation}\label{valeurbarA(a+eta)good}
  \bar A(\zeta+\eta)=A(\zeta+\eta), \, \bar A'(\zeta+\eta)=A'(\zeta+\eta).
\end{equation}
\end{lem}
\textbf{Proof of Lemma~\ref{modifBG}.}
Let us first consider the case where
\begin{equation}\label{memesigne}
  \left(B''(\zeta-\eta)+\frac{N-1}{\zeta-\eta}B'(\zeta-\eta)\right) \left(B''(\zeta+\eta)+\frac{N-1}{\zeta+\eta}B'(\zeta+\eta)\right)<0.
\end{equation}
Then, replacing if necessary $B$ by $-B$ and using \eqref{propSignG''}, we may assume
that
\begin{gather}
  \label{propertySignB''>a-eta}
  B''(z)+\frac{N-1}{z}B'(z)<0, \, \forall z \in [\zeta-\eta,\zeta),
  \\
\label{propertySignB''<a+eta}
  B''(z)+\frac{N-1}{z}B'(z)>0, \, \forall z \in (\zeta,\zeta+\eta).
\end{gather}
Let $\varphi \in C^\infty (-\infty,+\infty)$ be such that
\begin{gather}\label{propertyvarpi=1presde0}
  \varphi =1 \text{ in }[-1/2,1/2],
  \\
  \label{varphisupport}
  \varphi =0 \text{ in }(-\infty,-1]\cup[1,+\infty),
  \\
  \label{varphigeq0}
  \varphi(z) \in [0,1],\,  \forall z\in (-\infty,\infty).
\end{gather}
Let
\begin{equation}\label{defcalE}
  \mathcal{E}:=\left\{ \xi \in C^\infty([\zeta-\eta,\zeta+\eta]);\,
  \text{the support of $\xi$ is included in $(\zeta-\eta,\zeta+\eta)\setminus\{\zeta\}$}\right\}.
\end{equation}
The vector space $\mathcal{E}$ is equipped with the norm
\begin{equation}\label{defnormcalE}
  |\xi|:=\max\{|\xi(x)|;\, x\in [\zeta-\eta,\zeta+\eta]\}.
\end{equation}
For $\varepsilon\in \R$ and $\xi \in \mathcal{E}$, one
defines now $H_{\varepsilon,\xi}\in C^\infty([\zeta-\eta,\zeta+\eta])$  by,  if $\varepsilon\not = 0$,
\begin{equation}\label{defHepsilonnot0}
  H_{\varepsilon,\xi}(z):=\varepsilon^2(z-\zeta)^3\varphi\left(\frac{z-\zeta}{|\varepsilon|}\right)
  +\left(1-\varphi\left(\frac{z-\zeta}{|\varepsilon|}\right)\right)\left(
  B''(z)+\frac{N-1}{z}B'(z)+\xi(z)\right),
\end{equation}
for every $z\in [\zeta-\eta,\zeta+\eta]$ and
\begin{equation}\label{defHepsilon0}
  H_{0,\xi}(z):=
  B''(z)+\frac{N-1}{z}B'(z)+\xi(z), \, \forall z\in [\zeta-\eta,\zeta+\eta].
\end{equation}
We then define $\bar B:= B_{\varepsilon,\xi}\in C^\infty([\zeta-\eta,\zeta+\eta])$ by requiring
\begin{gather}\label{eqbarB}
  B_{\varepsilon,\xi}''(z)+\frac{N-1}{z}B_{\varepsilon,\xi}'(z)= H_{\varepsilon,\xi}(z),  \\
  \label{barBinit}
  B_{\varepsilon,\xi}(\zeta-\eta)=B(\zeta-\eta),\,
  B_{\varepsilon,\xi}'(\zeta-\eta)=B'(\zeta-\eta).
\end{gather}
Let $C_{\varepsilon,\xi}\in C^0([\zeta-\eta,\zeta+\eta])$ be defined by
\begin{equation}\label{defCepsilonzeta}
  C_{\varepsilon,\xi}(z):=-\left(B_{\varepsilon,\xi}''(z)+\frac{N-1}{z} B_{\varepsilon,\xi}'(z)\right)^{1/3}= -H_{\varepsilon,\xi}(z)^{1/3}.
\end{equation}
Note that by \eqref{propertyvarpi=1presde0}, \eqref{defHepsilonnot0}, and  \eqref{defCepsilonzeta}, if  $\varepsilon\ne0$,
\begin{equation}\label{valueCprime}
C'_{\varepsilon,\xi}(\xi)=-|\varepsilon|^{2/3}\not=0.
\end{equation}

Using \eqref{varphisupport}, \eqref{defcalE}, \eqref{defHepsilonnot0}, \eqref{defHepsilon0}, and  \eqref{eqbarB}, one sees that, if  $\varepsilon < \eta$ (which is assumed from now on),  $B_{\varepsilon,\xi}$ and $ B$ are both solutions to the second order differential equation
\begin{equation}\label{eqdiff-ordre-2}
Y''(z)+\frac{N-1}{z}Y'(z)= B''(z)+\frac{N-1}{z}B'(z)
\end{equation}
in a neighborhood of $\{\zeta-\eta,\zeta+\eta\}$ in $[\zeta-\eta,\zeta+\eta]$. In particular, by
\eqref{barBinit}, $B_{\varepsilon,\xi}$ and $B$ are equal in a neighborhood of $\zeta-\eta$ in $[\zeta-\eta,\zeta+\eta]$ and  \eqref{supportbarB} is equivalent to
\begin{equation}\label{barBgoonfin}
  B_{\varepsilon,\xi}(\zeta+\eta)=B(\zeta+\eta),\, B_{\varepsilon,\xi}'(\zeta+\eta)=B'(\zeta+\eta).
\end{equation}
Let $A_{\varepsilon,\xi} \in C^\infty ([\zeta-\eta,\zeta+\eta])$ be the solution of
\begin{gather}\label{equationbarepsilonalphaA}
  A_{\varepsilon,\xi}''+\frac{N-1}{z}A_{\varepsilon,\xi}' =-B_{\varepsilon,\xi}^3, \\
  \label{conditioninitbarepsilonalphaA}
  A_{\varepsilon,\xi}(\zeta-\eta)=A(\zeta-\eta),\, A_{\varepsilon,\xi}'(\zeta-\eta)=A'(\zeta-\eta).
\end{gather}
Let $\mathcal{F}:(-\eta,\eta)\times \mathcal{E}\rightarrow \R^4$ be defined
by
\begin{multline}\label{defcalF}
  \mathcal{F}(\varepsilon, \xi):=(B_{\varepsilon,\xi}(\zeta+\eta)-B(\zeta+\eta), B_{\varepsilon,\xi}'(\zeta+\eta)-B'(\zeta+\eta),
  \\
  A_{\varepsilon,\xi}(\zeta+\eta)-A(\zeta+\eta), A_{\varepsilon,\xi}'(\zeta+\eta)-A'(\zeta+\eta))\tr .
\end{multline}
One easily checks that
\begin{gather}\label{FC1}
  \mathcal{F} \text{ is of class }C^1,
  \\
  \label{F00}
  \mathcal{F}(0,0)=0.
\end{gather}

Let us assume, for the moment, that
\begin{equation}\label{F'onto}
  \frac{\partial \mathcal{F}}{\partial \xi}(0,0) \text{ is onto}.
\end{equation}
By \eqref{F'onto}, there exists a 4-dimensional subspace $\mathcal{E}_0$ of $\mathcal{E}$ such that
\begin{equation}\label{E0onto}
  \frac{\partial \mathcal{F}}{\partial \xi}(0,0)\mathcal{E}_0=\R^4.
\end{equation}
By \eqref{E0onto} and the implicit function theorem, there exists $\varepsilon_0\in (0,\eta)$ and a
 map $\xi :(-\varepsilon_0,\varepsilon_0)\rightarrow \mathcal{E}_0$ such that
\begin{gather}
\label{alpha0}
\xi(0)=0,
\\
\label{bonnefin}
  \mathcal{F}(\varepsilon, \xi(\varepsilon))=0, \,
  \forall \varepsilon \in (-\varepsilon_0,\varepsilon_0).
\end{gather}
From \eqref{propertySignB''>a-eta}, \eqref{propertySignB''<a+eta}, \eqref{defcalE}, \eqref{defnormcalE}, \eqref{defHepsilonnot0}, \eqref{defHepsilon0}, and \eqref{eqbarB}, one gets the existence
of $\varepsilon_1>0$ such that
\begin{gather}
\label{barB''bonsigne>}
B_{\varepsilon,\xi}''(z)+\frac{N-1}{z}B_{\varepsilon,\xi}'(z)<0,\,  \forall z \in [\zeta-\eta, \zeta), \, \forall \varepsilon \in [-\varepsilon_1,\varepsilon_1], \, \forall
\xi \in \mathcal{E}_0 \text{ such that } |\xi|\leqslant \varepsilon_1,
\\
\label{barB''bonsigne<}
B_{\varepsilon,\xi}''(z)+\frac{N-1}{z}B_{\varepsilon,\xi}'(z)>0,\, \forall z \in (\zeta,\zeta+\eta], \, \forall \varepsilon \in [-\varepsilon_1,\varepsilon_1], \, \forall
\xi \in \mathcal{E}_0 \text{ such that } |\xi|\leqslant \varepsilon_1.
\end{gather}
 From \eqref{propertyvarpi=1presde0}, \eqref{defHepsilonnot0}, and \eqref{eqbarB} one gets that, for every $\varepsilon \in (0,+\infty)$ and for every $\xi \in \mathcal{E}_0$, one has
\begin{equation}\label{presdea}
  B_{\varepsilon,\xi}''(z)+\frac{N-1}{z} B_{\varepsilon,\xi}'(z)=\varepsilon^2(z-\zeta)^3 \text{ if } |z-\zeta|\leqslant \varepsilon/2.
\end{equation}
 From \eqref{defCepsilonzeta},
\eqref{barB''bonsigne>}, \eqref{barB''bonsigne<}, and \eqref{presdea} one gets that,
for every $\varepsilon \in [-\varepsilon_1,\varepsilon_1]\setminus\{0\}$ and for every
$\xi \in \mathcal{E}_0$ such that  $|\xi|\leqslant \varepsilon_1$,
\begin{gather}\label{barb''1/3bon}
C_{\varepsilon,\xi} \in C^\infty([\zeta-\eta,\zeta+\eta]),
\\
\label{Cnotzerosaufzeta}
\left(C_{\varepsilon,\xi}(z)=0\right) \Leftrightarrow \left(z=\zeta\right),
\end{gather}
which, together with \eqref{E0onto} as above, \eqref{valueCprime}, \eqref{alpha0}, and \eqref{bonnefin}, conclude the proof of Lemma~\ref{modifBG} when \eqref{memesigne} holds.

It remains to prove \eqref{F'onto}. Simple computations show that
\begin{equation}\label{calculF'}
  \frac{\partial \mathcal{F}}{\partial \xi}(0,0)\xi = (x_1(\zeta+\eta),x_2(\zeta+\eta), x_3(\zeta+\eta), x_4(\zeta+\eta))\tr,
\end{equation}
where $x:[\zeta-\eta,\zeta+\eta]\rightarrow \R^4$ is the solution of
\begin{equation}\label{equationx}
  \dot x =K(t) x + \xi(t)e,
\end{equation}
with
\begin{equation}\label{defKdefb}
  K(t):=
  \begin{pmatrix}
  0&1&0&0
  \\
  0&-\frac{N-1}{t}&0&0
  \\
  0&0&0&1
  \\
  -3B^2(t)&0&0&-\frac{N-1}{t}
  \end{pmatrix}
  ,\,
  e:=
  \begin{pmatrix}
  0
  \\
  1
  \\
  0
  \\
  0
  \end{pmatrix}
  ,
\end{equation}
which satisfies
\begin{equation}\label{condinitx}
  x(\zeta-\eta)=0.
\end{equation}
Hence, using a standard density argument, \eqref{E0onto} comes from the following lemma.
\begin{lem}
\label{lemsys4control}
Let $\nu>0$, $\zeta>0$, $\eta >0$ be such that
$[\zeta-\eta,\zeta+\eta]\subset (0,+\infty)$.
Let $B\in C^\infty ([\zeta-\eta,\zeta+\eta])$ be such that
\begin{gather}\label{poureviternoncontrolable}
B\not \equiv 0.
\end{gather}
Then the control system \eqref{equationx}, where the state is $x\in \R^4$
and the control is $\xi \in \R$, is controllable on $[\zeta-\eta,\zeta+\eta]$, i.e. for every $X$ in $\R^4$ there exists $\xi \in L^\infty (\zeta-\eta,\zeta+\eta)$ such
that the solution of \eqref{equationx} and \eqref{condinitx} satisfies $x(\zeta+\eta)=X$.
\end{lem}
\textbf{Proof of Lemma~\ref{lemsys4control}.}
We use a classical result on the controllability of time-varying linear finite-dimensional
control systems (see e.g. \cite[Theorem 1.18]{2007-Coron-book}). One defines,
by induction on $i\in \mathbb{N}$, $e_i\in C^\infty([\zeta-\eta,\zeta+\eta])$ by requiring
\begin{gather}\label{defe0}
  e_0(t):= e, \, \forall t \in [\zeta-\eta,\zeta+\eta],\\
  e_i(t):=\dot e_{i-1}(t) -K(t) e_{i-1}(t), \, \forall t \in [\zeta-\eta,\zeta+\eta],
  \, \forall i \in \mathbb{N}\setminus\{0\}.
\end{gather}
Let $\theta\in C^\infty([\zeta-\eta,\zeta+\eta])$ be defined by
\begin{equation}\label{deftheta}
  \theta(t):= -\frac{N-1}{t},\, \forall t \in [\zeta-\eta,\zeta+\eta].
\end{equation}
Straightforward computations lead to
\begin{equation}\label{e1etc}
  e_1=
  \begin{pmatrix}
  -1
  \\
  -\theta
  \\
  0
  \\
  0
  \end{pmatrix}
  ,\,
  e_2=
  \begin{pmatrix}
  \theta
  \\
  -\dot \theta + \theta^2
  \\
  0
  \\
  -3B^2
  \end{pmatrix}
  ,\,
  e_3=
  \begin{pmatrix}
  2 \dot \theta - \theta^2
  \\
  -\ddot \theta +3\theta\dot \theta - \theta^3
  \\
  3B^2
  \\
  6B^2\theta -6 B\dot B
  \end{pmatrix}
  .
\end{equation}
From \eqref{defKdefb}, \eqref{defe0}, and \eqref{e1etc}, one gets
\begin{equation}\label{calculdet}
\text{det}(e_0,e_1,e_2,e_3)=9 B^4,
\end{equation}
which, with \eqref{poureviternoncontrolable} and \cite[Theorem 1.18]{2007-Coron-book},
concludes the proof of Lemma~\ref{lemsys4control}.
\cqfd

We now turn to the case where \eqref{memesigne} does not hold. Then, replacing if necessary $B$ by $-B$ and using \eqref{propSignG''}, we may assume
that
\begin{gather}
  \label{propertySignB''>0}
  B''(z)+\frac{N-1}{z}B'(z)>0, \, \forall z \in [\zeta-\eta,\zeta+\eta]\setminus\{\zeta\}.
\end{gather}
In the definition of $H_{\varepsilon,\xi}$ one replaces \eqref{defHepsilonnot0} by \begin{equation}\label{defHepsilonnew}
  H_{\varepsilon,\xi}(z):=\varepsilon^2\varphi\left(\frac{z-\zeta}{|\varepsilon|}\right)
  +\left(1-\varphi\left(\frac{z-\zeta}{|\varepsilon|}\right)\right)\left(
  B''(z)+\frac{N-1}{z}B'(z)+\xi(z)\right),
\end{equation}
and keeps \eqref{defHepsilon0}. Now \eqref{barB''bonsigne>} and \eqref{barB''bonsigne<} are  replaced by
\begin{gather}
\label{barC>0}
C_{\varepsilon,\xi}(z)>0,\,  \forall z \in [\zeta-\eta, \zeta+\eta], \, \forall \varepsilon \in [-\varepsilon_1,\varepsilon_1]\setminus\{0\}, \, \forall
\xi \in \mathcal{E}_0 \text{ such that } |\xi|\leqslant \varepsilon_1.
\end{gather}
Therefore, (compare with \eqref{Cnotzerosaufzeta}), provided that $\varepsilon\ne0$, one can see that $C_{\varepsilon,\xi}(z)\ne0$  for every $z \in [\zeta-\eta,\zeta+\eta]$ and consequently \eqref{C'enzerodeC} is satisfied. Moreover
\begin{gather}
\label{Csmooth}
C_{\varepsilon,\xi}
\in C^\infty([\zeta-\eta, \zeta+\eta]), \, \forall \varepsilon \in [-\varepsilon_1,\varepsilon_1]\setminus\{0\},
 \, \forall
\xi \in \mathcal{E}_0 \text{ such that } |\xi|\leqslant \varepsilon_1.
\end{gather}
which, together with \eqref{Csmooth}, \eqref{alpha0}, \eqref{bonnefin}, and \eqref{barC>0}, concludes the proof of Proposition~\ref{thmexistenceabc-stationary}.
\cqfd

\section{Proof of Theorem~\ref{thmexistenceabc} (time-varying case)}
\label{sectimevarying}
\setcounter{equation}{0}
In this section, we prove Theorem~\ref{thmexistenceabc}.
We define $\lambda \in C^{\infty}([-1,1])$ and $f_0 \in C^\infty([-1,1])$ by
\begin{equation}
\label{deflambda}
\lambda(t):=(1-t^2)^2,\, \forall t\in [-1,1],
\end{equation}
and
\begin{equation}
\label{deff0}
f_0(t):=\left\lbrace { \begin{tabular}{ll} \medskip
$e^{-\frac{1}{1-t^2}}$ & \text{ if }$|t|<1 $,\\ \medskip
$0$ & \text{ if }$t=0 $. \\ \end{tabular} } \right.
\end{equation}
Let $\varepsilon\in (0,1]$. For $r\in \R$ and $t\in (-1,1)$, we set
\begin{equation}\label{defz}
  z:=\frac{r}{\varepsilon \lambda(t)} \in [0,+\infty).
\end{equation}
Let $A$, $B$, and $C$ be as in Proposition~\ref{thmexistenceabc-stationary}. By \eqref{Cpair}, \eqref{valuebarCprochede0}, \eqref{Cpresde1sur2}, \eqref{barcpresde1}, and \eqref{mauvaisfinite}, there exist $p\in \mathbb{N}$ and $\rho_1$, $\rho_2$ \ldots $\rho_p$ in
$(-1,1)\setminus \{0\}$ such that
\begin{equation}\label{defrhoi}
  \left\{z\in (-1,1);\, C(z)=0\right\}=
  \left\{\rho_l;\, l\in\left\{1,2,\ldots,p\right\}\right\}.
\end{equation}
Let
\begin{equation}\label{defrho0}
  \rho_0:=\frac{1}{2}, \,\rho_{-1}:=-\frac{1}{2} .
\end{equation}
Let $\delta>0$ be such that
\begin{gather}\label{loinde0etde1}
  [\rho_l-\delta,\rho_l+\delta]\subset (-1,1)\setminus \{0\},\,  \forall l  \in\{-1,0,1,\ldots,p\},  \\
  \label{loinentreeux}
  [\rho_l-\delta,\rho_l+\delta]\cap [\rho_{l'}-\delta,\rho_{l'}+\delta]=\emptyset , \,\forall
  (l,l') \in \{-1,0,1,\ldots,p\}^2 \text{ such that  $l\not =l'$.}
\end{gather}

Let $\mathbb{D}:=\{(t,r)\in (-1,1)\times\R;\, |r|<\varepsilon\lambda(t)\}$. We look for $a:(t,r)\in \mathbb{D}\mapsto a(t,r)\in \R$ in the following form

\begin{gather}\label{shapeainside}
  a(t,r)=f_0(t)A(z)+\sum_{l=-1}^p\sum_{i=1}^{3}f_{il}(t)g_{il}(z),
\end{gather}
where the functions $f_{il}$, $g_{il}$ are to be determined with the requirement that
\begin{equation}\label{supportgil}
  \text{the support of $g_{il}$ is included in $(\rho_l-\delta,\rho_l+\delta)$,} \,
  \forall i\in \{1,2,3\},\, \forall l\in\{-1,0,1,\ldots,p\}.
\end{equation}
Then $b:(t,r)\in\mathbb{D}\mapsto b(t,r)\in \R$ is defined by
\begin{equation}\label{defb}
 b:=\left(a_t-a_{rr}-\frac{N-1}{r}a_{r}\right)^{1/3},
\end{equation}
and, on every open subset of $\mathbb{D}$ on which $b$ is of class $C^2$ and $b_{r}/r$ is bounded, $c$ is defined by
\begin{equation}\label{defc}
c:=\left(b_t-b_{rr}-\frac{N-1}{r}b_{r}\right)^{1/3}.
\end{equation}
 For $l\in \{-1,0,1,\ldots,p\}$, let $\Sigma_l\subset \R\times\R$ be defined by
\begin{equation}\label{xzonel}
\Sigma_l:=\left\{(t,r)\in (-1,1)\times \R; \, z\in (\rho_l-\delta,\rho_l+\delta)\right\}.
\end{equation}
Let us first study the case where, for some
\begin{equation}\label{propertybarl}
  \bar l \in \left\{1,2,\ldots,p\right\},
\end{equation}
$(t,r)\in \Sigma_{\bar l}$. By symmetry, we may only study  the case where $\rho_{\bar l}>0$.  Note that \eqref{propertybarl}, together with \eqref{C(1/2)notzero} and \eqref{defrhoi},
 implies that
\begin{equation}\label{rhonot=1/2}
\rho_{\bar l}\not =\frac{1}{2}.
\end{equation}
 From \eqref{loinentreeux}, \eqref{shapeainside},
\eqref{supportgil}, and \eqref{xzonel}, we have
\begin{equation}\label{valueazonebarl}
  a(t,r)=f_0(t)A(z)+\sum_{i=1}^{3}f_{i\bar l}(t)g_{i\bar l}(z).
\end{equation}
In order to simplify the notations,  we omit the index $\bar l$, and define $g_0$ by
\begin{equation}\label{defg0}
  g_0:=A.
\end{equation}
(This  definition is used all throughout this section.) Then, \eqref{valueazonebarl} now reads
\begin{equation}\label{valueazonebarlbarlomitted}
  a(t,r)=\sum_{i=0}^{3}f_{i}(t)g_{i}(z).
\end{equation}
Note that \eqref{eqBnulseulementen1/2}, \eqref{rhonot=1/2}, and \eqref{defg0} imply that
\begin{gather}
\label{Bnonnulenrho}
B(\rho)\not =0.
\end{gather}
Moreover, by \eqref{eqBC3}, \eqref{C'enzerodeC}, \eqref{defrhoi}, \eqref{propertybarl}, and  \eqref{defg0},
\begin{gather}\label{valeurnulleordre0}
\left(B^{(2)}+\frac{N-1}{z}B^{(1)}\right)(\rho)=0,
\\
\label{premierederiveenulle}
\left(B^{(2)}+\frac{N-1}{z}B^{(1)}\right)_z(\rho)=0,
\\
\label{deuxiemederiveenulle}
\left(B^{(2)}+\frac{N-1}{z}B^{(1)}\right)_{zz}(\rho)=0,
\\
\label{troisiemederiveenonnulle}
\left(B^{(2)}+\frac{N-1}{z}B^{(1)}\right)_{zzz}(\rho)\not=0.
\end{gather}
 To simplify the notations we assume that, for example,
\begin{gather}
\label{Bnegativeatrho}
B(\rho)<0,
\\
\label{troisiemederivee>0}
\left(B^{(2)}+\frac{N-1}{z}B^{(1)}\right)_{zzz}(\rho)<0.
\end{gather}
From \eqref{premierederiveenulle}, \eqref{deuxiemederiveenulle},
 \eqref{Bnegativeatrho}, and \eqref{troisiemederivee>0}, if $\delta \in (0,\rho)$ is small enough, there exists $\mu>0$ such that
\begin{gather}
\label{Bnonnulenvoisinagerho}
B(z)\leqslant -\mu , \, \forall z\in [\rho-\delta, \rho+\delta],
\\
\label{troisiemederivee>mu}
\left(B^{(2)}+\frac{N-1}{z}B^{(1)}\right)_{zzz}(z)\leqslant -\mu, \, \forall z\in [\rho-\delta, \rho+\delta].
\end{gather}
We now fix such a $\delta$.

 From \eqref{defb} and \eqref{valueazonebarlbarlomitted},
\begin{equation}\label{expressionb-1}
  b=
-\frac{1}{\varepsilon^{2/3} \lambda^{2/3}}\left(\sum_{i=0}^{3} \Big(f_ig_i^{(2)}
+\frac{N-1}{z}f_ig^{(1)}_i
+z\varepsilon^2\lambda\dot{\lambda}f_ig_i^{(1)}-\varepsilon^2\lambda^2\dot{f_i}g_i \Big)\right)^{1/3}.
\end{equation}
Let us denote by $M:\R\times \R^*\to\R$, $(t,z)\mapsto M(t,z)\in \R$, the function defined by:
\begin{multline}\label{defM}
  M(t,z):=\sum_{i=0}^{3} \Big(f_i(t)g_i^{(2)}(z)+\frac{N-1}{z}f_i(t)g^{(1)}_i(z)
  \\ +z\varepsilon^2\lambda(t)\dot{\lambda}(t)f_i(t)g_i^{(1)}(z)
  -\varepsilon^2\lambda^2(t)\dot{f_i}(t)g_i(z)\Big).
\end{multline}
For the moment, let us assume that
\begin{equation}\label{Mnot0}
  M(t,z)\not =0, \, \forall (t,z)\in (-1,1)\times (\rho-\delta, \rho +\delta).
\end{equation}
Using \eqref{defz}, \eqref{defc},
\eqref{expressionb-1}, \eqref{defM}, and straightforward computations, one gets, on the open set of the $(t,r)\in \Sigma$ such that $M(t,z)\not=0$,
\begin{equation}
\label{valuec}
9\varepsilon^{8/3}\lambda^{8/3}c^3=\nu,
\end{equation}
with
\begin{equation}\label{defnu}
\nu:=\frac{1}{M^{2/3}}\left(
 3 M_{zz}
-2 \frac{M_{z}^2}{M}
+\frac{3(N-1)}{ z} M_z +6 \varepsilon^2\lambda \dot \lambda M
- 3\varepsilon^2\lambda^2 M_t+ 3z \varepsilon^2 \lambda \dot \lambda M_z\right).
\end{equation}

The idea is to construct the $f_i$'s and the $g_{i}$'s in order to have a precise knowledge of the places where $\nu$ vanishes and the order of the vanishing so that $\nu$ is the cube of a $C^\infty$ function. More precisely, we are are going to check that one can construct the $f_i$'s and the $g_{i}$'s so that, at least if $\varepsilon\in (0,1]$ is small enough,
\begin{gather}
\label{nurho=0}
\nu(t,\rho)=0, \, \forall t \in (-1,1),
\\
\label{nuzrho=0}
\nu_z(t,\rho)=0, \, \forall t \in (-1,1),
\\
\label{nuzzrho=0}
\nu_{zz}(t,\rho)=0, \, \forall t \in (-1,1),
\\
\label{nuzzzrhonot =0}
\nu_{zzz}(t,\rho)>0, \, \forall t \in (-1,1).
\end{gather}

 From \eqref{defM}, one has
\begin{equation}
\label{valueMz}
M_z=\sum_{i=0}^{3}\Big( f_ig_i^{(3)}+\varepsilon^2(\dot{\lambda}f_i-\lambda\dot{f_i})\lambda g_i^{(1)}+\varepsilon^2 z\lambda \dot{\lambda}f_ig_i^{(2)}+\frac{N-1}{z}f_ig_i^{(2)}-\frac{N-1}{z^2}f_ig_i^{(1)}\Big ),
\end{equation}
\begin{equation}
\label{valueMzz}
\begin{array}{rcl}
M_{zz}&=&\sum_{i=0}^{3}\Big( f_ig_i^{(4)}+\varepsilon^2(2\dot{\lambda}f_i-\lambda\dot{f_i})\lambda g_i^{(2)}+\varepsilon^2 z\lambda \dot{\lambda}f_ig_i^{(3)}
\\
&&+\frac{N-1}{z}f_ig_i^{(3)}-\frac{2(N-1)}{z^2}f_ig_i^{(2)}+\frac{2(N-1)}{z^3}f_ig_i^{(1)} \Big).
\end{array}
\end{equation}
We impose that
\begin{equation}
\label{propertygiderineq4}
g_i^{(j)}(\rho)=
\left\{
\begin{array}{rl}
1 & \text{ if $i=1$ and $j=4$},\\
0 & \text{ if $1\leq i\leq 3$, $0\leq j\leq 4$ and $(i,j)\not = (1,4)$.} \end{array}
\right.
\end{equation}
 From \eqref{defM}, \eqref{valueMz}, \eqref{valueMzz}, and \eqref{propertygiderineq4}, we have
\begin{gather}
\label{expMatrho}
M(\cdot ,\rho)=f_0g_0^{(2)}(\rho)+\frac{N-1}{\rho}f_0g_0^{(1)}(\rho) +\varepsilon^2 \rho\lambda\dot{\lambda}f_0g_0^{(1)}(\rho)
-\varepsilon^2 \lambda^2\dot{f_0}g_0(\rho),
\end{gather}
\begin{equation}
\label{expMzatrho}
\begin{array}{rcl}
M_z(\cdot,\rho)&=&f_0g_0^{(3)}(\rho)\displaystyle +\frac{N-1}{\rho}f_0g_0^{(2)}(\rho)-\frac{N-1}{\rho^2}f_0g_0^{(1)}(\rho)
\\
&&
\displaystyle +\varepsilon^2(\dot{\lambda}f_0-\lambda\dot f_0)\lambda g_0^{(1)}(\rho)+\varepsilon^2\rho\lambda \dot{\lambda} f_0g_0^{(2)}(\rho)
,
\end{array}
\end{equation}
\begin{gather}
\label{expMzzatrho}
\begin{array}{rcl}
M_{zz}(\cdot,\rho)&=&f_0g_0^{(4)}(\rho)\displaystyle
+\frac{N-1}{\rho}f_0g_0^{(3)}(\rho) -\frac{2(N-1)}{\rho^2}f_0g_0^{(2)}(\rho)+\frac{2(N-1)}{\rho^3}f_0g_0^{(1)}(\rho)
\\
&&
+f_1+\varepsilon^2
(2\dot{\lambda}f_0-\lambda\dot{f_0})\lambda g_0^{(2)}(\rho)+\varepsilon^2 \rho\lambda \dot{\lambda} f_0g_0^{(3)}(\rho).
\end{array}
\end{gather}
  From \eqref{eqAB3}, \eqref{deflambda}, \eqref{deff0}, \eqref{defg0}, \eqref{Bnegativeatrho}, \eqref{expMatrho}, and \eqref{expMzatrho}, one has,
at least if  $\varepsilon >0$ is small enough, which is from now on assumed,
\begin{equation}
\forall t \in(-1,1), M(t,\rho)>  0.
\end{equation}
Then, for  $z=\rho$, one has
\begin{equation} \label{valuedenuarho}
\begin{array}{rcl}
\displaystyle
\nu(.,\rho)&=\frac{1}{M^{2/3}(.,\rho)}&\Big(3M_{zz}(.,\rho)-2\frac{M_z^2(.,\rho)}{M(.,\rho)}+\frac{3(N-1)}{\rho}M_z(.,\rho)
+6\varepsilon ^2\lambda\dot{\lambda}M_z(.,\rho)
\\
&&\displaystyle
-
3\varepsilon ^2\lambda^2M_t(.,\rho)+3\rho\varepsilon^2\lambda\dot{\lambda}M_z(.,\rho)\Big).
\end{array}
\end{equation}
We then choose to define $f_1: t\in (-1,1) \mapsto f_1(t)\in \R$  by
\begin{equation} \label{deff1}
\begin{array}{rcl}\displaystyle
f_1&:=&-f_0g_0^{(4)}(\rho)-\frac{N-1}{\rho}f_0g_0^{(3)}(\rho)+\frac{2(N-1)}{\rho^2}f_0g_0^{(2)}(\rho)
-\frac{2(N-1)}{\rho^3}f_0g_0^{(1)}(\rho)
\\
&&\displaystyle
-\varepsilon ^2(2\dot{\lambda}f_0-\lambda\dot{f_0})\lambda g_0^{(2)}(\rho)-\varepsilon ^2\rho\lambda \dot{\lambda} f_0g_0^{(3)}(\rho)
\\
&&\displaystyle
+\frac{1}{3}\Big(2\frac{M_z^2(.,\rho)}{M(.,\rho)}-\frac{3(N-1)}{\rho}M_z(.,\rho)
-6\varepsilon ^2\lambda\dot{\lambda}M(.,\rho)
\\
&&\displaystyle
+3\varepsilon ^2\lambda^2M_t(.,\rho)
-3\rho\varepsilon ^2\lambda\dot{\lambda}M_z(.,\rho)\Big).
 \end{array}
\end{equation}
Note that, even if $M$ depends on $f_1$, $f_2$, and $f_3$, the right hand side of \eqref{deff1} does not depend on $f_1$,
$f_2$, and $f_3$, and
$f_1$ is indeed well-defined by \eqref{deff1}. This definition of $f_1$, together with \eqref{expMzzatrho} and \eqref{valuedenuarho}, implies that \eqref{nurho=0} holds. (In fact, $f_1$ is defined by \eqref{deff1} precisely
in order to have \eqref{nurho=0}.) From \eqref{deflambda}, \eqref{deff0}, \eqref{valeurnulleordre0}, \eqref{Bnegativeatrho}, \eqref{expMatrho}, \eqref{expMzatrho}, and \eqref{deff1},
we obtain the existence of two polynomials $p_1(\varepsilon^2,t)$ and $q_1(\varepsilon^2,t)$ in the variables $\varepsilon^2$ and $t$
 such that
\begin{gather}
\label{f1=epsilonq1}
f_1(t)=\varepsilon ^2 \frac{p_1(\varepsilon^2,t)}{1+ q_1(\varepsilon^2,t)}f_0(t) ,\, \forall t \in (-1,1).
\end{gather}

In order to simplify the notations, we  set:
\begin{equation} \label{defK}
 \begin{array}{rcl}\displaystyle
K(t,z)&:=&-\frac{2M_z(t,z)^2}{M(t,z)}
+\frac{3(N-1)}{z}M_z(t,z)+6\varepsilon^2\lambda\dot{\lambda}M(t,z)
\\[2mm]
&&
\displaystyle-3\varepsilon^2\lambda^2M_t(t,z)+3z\varepsilon^2\lambda\dot{\lambda}M_z(t,z).
 \end{array}
\end{equation}
We then have

\begin{equation} \label{r1}
 \begin{array}{l}\displaystyle
\nu=M^{-\frac{2}{3}}(3M_{zz}+K).
 \end{array}
\end{equation}
Differentiating this equality with respect to $z$, we obtain

\begin{equation} \label{nuz=}
 \begin{array}{l}\displaystyle
\nu_z=M^{-\frac{5}{3}}(3MM_{zzz}+MK_z-2M_zM_{zz}-\frac{2}{3}M_zK).
 \end{array}
\end{equation}
Differentiating \eqref{defK} with respect to $z$, we get
\begin{equation} \label{valueKz}
\begin{array}{rcl}
K_z&=&-\frac{4M_zM_{zz}}{M^2}+\frac{2M_z^3}{M^3}
+\frac{3(N-1)}{z}M_{zz}-\frac{3(N-1)}{z^2}M_{z}
\\[2mm]
&&
\displaystyle +9\varepsilon^2\lambda\dot{\lambda}M_z-3\varepsilon^2\lambda^2M_{tz}
+3\varepsilon ^2z\lambda\dot{\lambda}M_{zz}.
\end{array}
\end{equation}

Then, differentiating \eqref{valueMzz} with respect to $z$, we have
\begin{equation} \label{valueMzzz}
\begin{array}{rcl}
\displaystyle
M_{zzz}&=\sum_{i=0}^3&\Big(f_ig_i^{(5)}+\frac{N-1}{z}f_ig_i^{(4)}
-\frac{3(N-1)}{z^2}f_ig_i^{(3)} +\frac{6(N-1)}{z^3}f_ig_i^{(2)}
\\
&&
\displaystyle
-\frac{6(N-1)}{z^4}f_ig_i^{(1)}
+ \varepsilon^2z\lambda \dot{\lambda}f_ig_i^{(4)}+\varepsilon^2(3\dot{\lambda}f_i-\lambda \dot f_i)\lambda g_i^{(3)}\Big).
 \end{array}
\end{equation}
We impose that
\begin{equation}
\label{defgi5}
g_i^{(5)}(\rho)= \left\lbrace { \begin{tabular}{ll} \medskip
$1$ & \text{ if }$i=2$,\\
$0$ & \text{ if $i\in \{1,3\}$}. \end{tabular} } \right.
\end{equation}
 From \eqref{propertygiderineq4}, \eqref{valueMzzz}, and \eqref{defgi5}, we have
\begin{equation} \label{expMzzzatrho}
\begin{array}{rcl}
\displaystyle
M_{zzz}(.,\rho)&=&f_0g_0^{(5)}(\rho)+\frac{N-1}{\rho}f_0g_0^{(4)}(\rho)
-\frac{3(N-1)}{\rho^2}f_0g_0^{(3)}(\rho)+\frac{6(N-1)}{\rho^3}f_0g_0^{(2)}(\rho)
\\
&&
\displaystyle
-\frac{6(N-1)}{\rho^4}f_0g_0^{(1)}(\rho)+f_2+\frac{N-1}{\rho}f_1
\\
&&+\varepsilon^2(3\dot{\lambda}f_0-\lambda\dot{f_0})\lambda g_0^{(3)}(\rho)
+
\varepsilon^2\rho\lambda \dot{\lambda}f_0g_0^{(4)}(\rho)
+\varepsilon ^2\rho\lambda \dot{\lambda}f_1.
 \end{array}
\end{equation}
We then define $f_2: t\in (-1,1) \mapsto f_2(t)\in \R$  by
\begin{equation} \label{deff2}
 \begin{array}{rcl}
\displaystyle
f_2&:=&-f_0g_0^{(5)}(\rho)
\displaystyle
-\frac{N-1}{\rho}f_0g_0^{(4)}(\rho) +\frac{3(N-1)}{\rho^2}f_0g_0^{(3)}(\rho)
\\
&&
\displaystyle
-\frac{6(N-1)}{\rho^3}f_0g_0^{(2)}(\rho)+\frac{6(N-1)}{\rho^4}f_0g_0^{(1)}(\rho)
-\frac{N-1}{\rho}f_1 \\
&&
\displaystyle
+\frac{1}{3M(.,\rho)}(-M(.,\rho)K_z(.,\rho)+2M_z(.,\rho)M_{zz}(.,\rho)+\frac{2}{3}M_z(.,\rho)K(.,\rho))
\\
&&
-\varepsilon^2\rho\lambda \dot{\lambda}f_0g_0^{(4)}(\rho)-\varepsilon^2\rho\lambda \dot{\lambda}f_1
-\varepsilon ^2(3\dot{\lambda}f_0-\lambda\dot{f_0})\lambda g_0^{(3)}(\rho).
\end{array}
\end{equation}
Note that, again, even if $M$ depends on $f_2$ and $f_3$, the right hand side of \eqref{deff2} does not depend on $f_2$ and $f_3$ (it depends on $f_1$, however $f_1$ is already defined in \eqref{deff1}), and
$f_2$ is indeed well defined by \eqref{deff2}. This definition of $f_2$, together with \eqref{nuz=} and \eqref{expMzzzatrho}, implies \eqref{nuzrho=0}.  From \eqref{deflambda}, \eqref{deff0}, \eqref{premierederiveenulle},
 \eqref{expMatrho}, \eqref{expMzatrho}, \eqref{expMzzatrho}, \eqref{f1=epsilonq1},
  \eqref{defK}, \eqref{valueKz}, and \eqref{deff2},
we obtain the existence of two polynomials $p_2(\varepsilon^2,t)$ and $q_2(\varepsilon^2,t)$ in the variables $\varepsilon^2$ and $t$
such that
\begin{gather}
\label{f2=epsilonq2}
f_2(t)=\varepsilon^2 \frac{p_2(\varepsilon^2,t)}{1+ q_2(\varepsilon^2,t)}f_0(t) ,\, \forall t \in (-1,1).
\end{gather}
Differentiating \eqref{nuz=} with respect to $z$, we obtain

\begin{equation} \label{expressionnuzz}
 \begin{array}{rcl}
\nu_{zz}&=&
\displaystyle M^{-\frac{8}{3}}\Big(-4MM_zM_{zzz}-\frac{7}{3}MM_zK_
z+\frac{10}{3}M_z^2M_{zz}+\frac{10}{9}M_z^2K
\\
&& \displaystyle
+3M^2M_{zzzz}
+MM_zK_z+M^2K_{zz}
-2MM_{zz}^2-\frac{2}{3}MM_{zz}K\Big).
 \end{array}
\end{equation}
Differentiating \eqref{valueKz} with respect to $z$, we obtain
\begin{equation}
\label{Kzz=}
 \begin{array}{rcl}
\displaystyle
K_{zz}&=&-\frac{4M_{zz}^2}{M}-\frac{4M_zM_{zzz}}{M}+\frac{10M_z^2M_{zz}}{M^2}-\frac{4M_z^4}{M^3}+\frac{6(N-1)}{z^3}M_z
\\[4mm]
&& \displaystyle
-\frac{6(N-1)}{z^2}M_{zz}+\frac{3(N-1)}{z}M_{zzz}
+12\varepsilon^2\lambda\dot{\lambda}M_{zz}-3\varepsilon^2\lambda^2M_{tzz}+3z\varepsilon^2\lambda\dot{\lambda}M_{zzz}.
 \end{array}
\end{equation}
Differentiating \eqref{valueMzzz} with respect to $z$, one has

\begin{equation} \label{Mzzzz=}
 \begin{array}{rll}\displaystyle
M_{zzzz}&=\sum_{i=0}^3 \Big(&f_ig_i^{(6)}+\frac{N-1}{z}f_ig_i^{(5)}
-\frac{4(N-1)}{z^2}f_ig_i^{(4)}
\\
[2mm]
&&\displaystyle
+ \frac{12(N-1)}{z^3}f_ig_i^{(3)}-\frac{24(N-1)}{z^4}f_ig_i^{(2)}
+\frac{24(N-1)}{z^5}f_ig_i^{(1)}
\\
[2mm]
&&\displaystyle
+\varepsilon^2(4\dot{\lambda}f_i-\lambda\dot{f_i})\lambda g_i^{(4)}+\varepsilon^2 z\lambda\dot{\lambda} f_ig_i^{(5)}
\Big).
\end{array}
\end{equation}
We then impose
\begin{equation}g_i^{(6)}(\rho)= \left\lbrace { \begin{tabular}{ll}
$1$ & \text{ if }$i=3$,\\
$0$ &  \text{ if $i\in \{1,2\}$}. \\ \end{tabular} } \right.
\end{equation}
Evaluating $M_{zzzz}$ at $z=\rho$ in \eqref{Mzzzz=} gives
\begin{equation} \label{valueMzzzzatrho}
 \begin{array}{rcl}\displaystyle
M_{zzzz}(.,\rho)&=&f_0g_0^{(6)}(\rho)+f_3+\frac{N-1}{\rho}f_0g_0^{(5)}(\rho) -\frac{4(N-1)}{\rho^2}f_0g_0^{(4)}(\rho)
\\
\displaystyle
&&+\frac{12(N-1)}{\rho^3}f_0g_0^{(3)}(\rho)
-\frac{24(N-1)}{\rho^4}f_0g_0^{(2)}(\rho)+\frac{24(N-1)}{\rho^5}f_0g_0^{(1)}(\rho)
\\ \noalign{\medskip}\displaystyle
&&
\displaystyle
+\frac{N-1}{\rho}f_2-\frac{4(N-1)}{\rho^2}f_1
+\varepsilon^2(4\dot{\lambda}f_0-\lambda\dot{f_0})\lambda g_0^{(4)}(\rho)
\\
&&+\varepsilon^2(4\dot{\lambda}f_1-\lambda\dot{f_1})\lambda+\varepsilon^2 \rho\lambda\dot{\lambda} f_0g_0^{(5)}(\rho)+\varepsilon^2\rho\lambda\dot{\lambda} f_2.
\end{array}
\end{equation}
Then, we define $f_3: t\in (-1,1) \mapsto f_3(t)\in \R$ by
\begin{equation} \label{deff3}
 \begin{array}{rcl}\displaystyle
f_3&:=&-f_0g_0^{(6)}(\rho)-\frac{N-1}{\rho}f_0g_0^{(5)}(\rho)
+\frac{4(N-1)}{\rho^2}f_0g_0^{(4)}(\rho) -\frac{12(N-1)}{\rho^3}f_0g_0^{(3)}(\rho)
\\
&&
\displaystyle
+\frac{24(N-1)}{\rho^4}f_0g_0^{(2)}(\rho) -\frac{24(N-1)}{\rho^5}f_0g_0^{(1)}(\rho) -\frac{N-1}{\rho}f_2  +\frac{4(N-1)}{\rho^2}f_1
\\
&&
\displaystyle
-\varepsilon^2 \rho\lambda\dot{\lambda} f_0g_0^{(5)}(\rho) -\varepsilon^2(4\dot{\lambda}f_0-\lambda\dot{f_0})\lambda g_0^{(4)}(\rho)-\varepsilon^2\rho\lambda\dot{\lambda} f_2-\varepsilon^2(4\dot{\lambda}f_1-\lambda\dot{f_1})\lambda
\\
&&\displaystyle
+\frac{1}{3M^2(.,\rho)}\Big(4M(.,\rho)M_z(.,\rho)M_{zzz}(.,\rho)
+\frac{7}{3}M(.,\rho)M_z(.,\rho)K_z(.,\rho)
\\
&&\displaystyle
-\frac{10}{3}M_z^2(.,\rho)M_{zz}(.,\rho)-\frac{10}{9}M_z^2(.,\rho)K(.,\rho)
-M(.,\rho)M_z(.,\rho)K_z(.,\rho)
\\
&& \displaystyle
-M^2(.,\rho)K_{zz}(.,\rho)
+2M(.,\rho)M_{zz}^2(.,\rho)
+\frac{2}{3}M(.,\rho)M_{zz}(.,\rho)K(.,\rho)\Big).
\end{array}
\end{equation}
Once more,  even if $M$ depends on $f_3$, the right hand side of \eqref{deff3} does not depend on $f_3$, and
$f_3$ is indeed well defined by \eqref{deff3}. This definition of $f_3$, together with \eqref{expressionnuzz}
and \eqref{valueMzzzzatrho}, implies that \eqref{nuzzrho=0} holds.
 From \eqref{deflambda}, \eqref{deff0}, \eqref{deuxiemederiveenulle}, \eqref{expMatrho},
 \eqref{expMzatrho},  \eqref{expMzzatrho}, \eqref{f1=epsilonq1}, \eqref{defK},
 \eqref{valueKz}, \eqref{expMzzzatrho}, \eqref{f2=epsilonq2},
 \eqref{Kzz=},  and \eqref{deff3},
we obtain the existence of two polynomials $p_3(\varepsilon^2,t)$ and $q_3(\varepsilon^2,t)$ in the variables $\varepsilon^2$ and $t$,
 such that
\begin{gather}
\label{f3=epsilonq3}
f_3(t)=\varepsilon^2 \frac{p_3(\varepsilon^2,t)}{1+ q_3(\varepsilon^2,t)}f_0(t) ,\, \forall t \in (-1,1).
\end{gather}

We are now in a position to analyse the regularity of $a$,  $b$, and $c$ on $\Sigma$.
Let us first point out that, by \eqref{deflambda}, \eqref{deff0}, \eqref{shapeainside}, \eqref{f1=epsilonq1}, \eqref{f2=epsilonq2}, and  \eqref{f3=epsilonq3}, there exists $\psi^a :
(t,z)\in [-1,1]\times [\rho-\delta,\rho+\delta]\mapsto \psi^a(t,z)\in \R $
of class $C^\infty$  such that
\begin{equation}\label{areguliersurSigma-cas1}
 a(t,r)=f_0(t)\psi^a(t,z), \forall (t,r)\in \Sigma.
\end{equation}
In particular, $a$ is of class $C^\infty$ in $\Sigma$.
 From \eqref{eqAB3}, \eqref{deflambda}, \eqref{deff0}, \eqref{defg0}, \eqref{expressionb-1},  \eqref{Bnonnulenvoisinagerho}, \eqref{f1=epsilonq1},  \eqref{f2=epsilonq2}, and \eqref{f3=epsilonq3}, we get that, at least if $\varepsilon>0$ is small enough, there there exists
 $\psi^b :
(t,z)\in [-1,1]\times [\rho-\delta,\rho+\delta]\mapsto \psi^b(t,z)\in \R $
of class  $ C^\infty$  such that
\begin{gather}
b<0 \text{ in }\Sigma,
\\
\label{breguliersurSigma-cas1}
b(t,r)=\lambda^{-2/3}f_0(t)^{1/3}\psi^b(t,z),
\forall (t,r)\in \Sigma.
\end{gather}
In particular, $b$ is of class $C^\infty$ in $\Sigma$.

Let us now study $c$. Differentiating \eqref{expressionnuzz} with respect to $z$ one gets
\begin{equation} \label{expressionnuzzz}
 \begin{array}{rcl}
\displaystyle
\nu_{zzz}&=M^{-\frac{11}{3}}\Big(&10MM_z^2M_{zzz}+\frac{10}{3}MM_z^2K_z-\frac{80}{9}M_z^3M_{zz}
-\frac{80}{27}M_z^3K-6M^2M_zM_{zzzz}
\\ \displaystyle
&& -2M^2M_zK_{zz} +10MM_zM_{zz}^2+\frac{10}{3}MM_zM_{zz}K-8M^2M_{zz}M_{zzz}
\\ \displaystyle
&& -2M^2M_{zz}K_z+3M^3M_{zzzzz}+M^3K_{zzz} -\frac{2}{3}M^2M_{zzz}K\Big).
 \end{array}
\end{equation}
Differentiating \eqref{Mzzzz=} with respect to $z$, we get
\begin{equation} \label{Mzzzzz=}
 \begin{array}{rll}\displaystyle
M_{zzzzz}&=\sum_{i=0}^3 \Big(&f_ig_i^{(7)}+\frac{N-1}{z}f_ig_i^{(6)}
-\frac{5(N-1)}{z^2}f_ig_i^{(5)} + \frac{20(N-1)}{z^3}f_ig_i^{(4)}
\\
[2mm]
&&\displaystyle
-\frac{60(N-1)}{z^4}f_ig_i^{(3)}
+\frac{120(N-1)}{z^5}f_ig_i^{(2)}-\frac{120(N-1)}{z^6}f_ig_i^{(1)}
\\
[2mm]
&&\displaystyle
+\varepsilon^2(5\dot{\lambda}f_i-\lambda\dot{f_i})\lambda g_i^{(5)}+\varepsilon^2 z\lambda\dot{\lambda} f_ig_i^{(6)}
\Big).
\end{array}
\end{equation}
Differentiating \eqref{Kzz=} with respect to $z$, we get
\begin{equation}
\label{Kzzz=}
 \begin{array}{rcl}
\displaystyle
K_{zzz}&=&\frac{6M_{z}^2M_{zz}}{M}-\frac{4M_zM_{zzzz}}{M}+\frac{24M_zM_{zz}^2}{M^2}
+\frac{4M_{zz}M_{zzz}}{M}-\frac{36M_z^3M_{zz}}{M^3}+\frac{12M_{z}^5}{M^4}
\\[4mm]
&& \displaystyle
-\frac{18(N-1)}{z^4}M_z +\frac{18(N-1)}{z^3}M_{zz}-\frac{9(N-1)}{z^2}M_{zzz}+\frac{3(N-1)}{z}M_{zzzz}
\\[4mm]
&& \displaystyle
+15\varepsilon^2\lambda\dot{\lambda}M_{zzz}-3\varepsilon^2\lambda^2M_{tzzz}+3z\varepsilon^2\lambda\dot{\lambda}M_{zzzz}.
 \end{array}
\end{equation}
 From  \eqref{deflambda}, \eqref{deff0}, \eqref{troisiemederivee>mu}, \eqref{defM}, \eqref{valuec},
 \eqref{defnu}, \eqref{nurho=0}, \eqref{nuzrho=0}, \eqref{nuzzrho=0}, \eqref{valueMz},
 \eqref{valueMzz}, \eqref{f1=epsilonq1},
 \eqref{defK},  \eqref{valueKz}, \eqref{valueMzzz},
   \eqref{f2=epsilonq2}, \eqref{Kzz=},  \eqref{f3=epsilonq3}, \eqref{expressionnuzzz}, \eqref{Mzzzzz=}, and \eqref{Kzzz=},
 one gets the existence of $\phi :
(t,z)\in [-1,1]\times [\rho-\delta,\rho+\delta]\mapsto \phi (t,z)\in \R $
of class  $ C^\infty$  such that
\begin{gather}
\label{breguliersurSigma-step1}
c^3(t,r)=\lambda^{8/3}f_0(t)^{1/3}\phi(t,z),
\forall (t,r)\in \Sigma,
\\
\label{phicderivee0}
\phi(t,\rho)=0,\, \forall t \in [-1,1],
\\
\label{phicderivee1}
\partial_z\phi(t,\rho)=0,\, \forall t \in [-1,1],
\\
\label{phicderivee2}
\partial^2_{zz}\phi(t,\rho)=0,\, \forall t \in [-1,1],
\\
\label{phicderivee3}
\partial^3_{zzz}\phi(t,z)>0,\, \forall (t,z)\in [-1,1]\times [\rho-\delta,\rho+\delta].
\end{gather}
Let $\tilde \phi :
(t,z)\in [-1,1]\times [\rho-\delta,\rho+\delta]\mapsto \tilde \phi (t,z)\in \R $  be defined by
\begin{gather}
\label{deftildephihorsrho}
\tilde \phi (t,z):=\frac{1}{2}\int_0^1(1-s)^2\partial^3_{zzz}\phi(t,\rho+s(z-\rho)))ds,\, \forall
 (t,z)\in [-1,1]\times [\rho-\delta,\rho+\delta].
\end{gather}
Then, $\tilde \phi$ is of class $C^\infty$ on $[-1,1]\times [\rho-\delta,\rho+\delta]$ and, using
\eqref{phicderivee0}, \eqref{phicderivee1}, \eqref{phicderivee2}, and \eqref{phicderivee3},
\begin{gather}
\label{phitildephi}
\phi(t,z)= (z-\rho)^3\tilde \phi (t,z), \, \forall (t,z)\in [-1,1]\times [\rho-\delta,\rho+\delta],
\\
\label{tildephi>0}
\tilde \phi (t,z)>0,\, \forall (t,z)\in [-1,1]\times [\rho-\delta,\rho+\delta].
\end{gather}
Let $\psi^c : (t,z)\in [-1,1]\times [\rho-\delta,\rho+\delta]\mapsto \psi^c(t,z)\in \R $ be defined by
\begin{equation}\label{defpsic}
\psi^c (t,z):= (z-\rho)\tilde \phi (t,z)^{1/3},\, \forall (t,z)\in [-1,1]\times [\rho-\delta,\rho+\delta].
\end{equation}
By \eqref{breguliersurSigma-step1}, \eqref{phitildephi}, \eqref{tildephi>0}, and \eqref{defpsic}, one gets that
\begin{gather}
\label{psicregular}
\psi^c \in C^\infty([-1,1]\times [\rho-\delta,\rho+\delta]),
\\
\label{creguliersurSigma-cas1}
c(t,r)=\lambda^{-8/9}f_0(t)^{1/9}\psi^c(t,z),
\forall (t,r)\in \Sigma.
\end{gather}
In particular, $c$ is of class $C^\infty$ in $\Sigma$.

Let us now study the case $l\in \{-1,0\}$, i.e. $\rho_l=1/2$ or $\rho_l=-1/2$. By symmetry, we may assume that $l=0$
so that $\rho_l=1/2$. This case is simpler than the previous one. It is already
 treated in \cite{2010-Coron-Guerrero-Rosier-SICON}, except that we now have  to take care of $c$.
 So, we will only briefly sketch the arguments. By \eqref{C(1/2)notzero} we may impose on $\delta$ to
be small enough so that
\begin{gather}
\label{C>0pourdeltassezpetit}
C(z)>0,\, \forall z\in [(1/2)-\delta,(1/2)+\delta].
\end{gather}
We now define (see \eqref{expressionb-1} and compare with \eqref{defnu})
\begin{equation}\label{newnu}
\nu := \sum_{i=0}^{3}\Big( f_ig_i^{(2)}
+\frac{N-1}{z}f_ig^{(1)}_i
+z\varepsilon^2\lambda\dot{\lambda}f_ig_i^{(1)}-\varepsilon^2\lambda^2\dot{f_i}g_i\Big).
\end{equation}
We still want to ensure that \eqref{nurho=0} to \eqref{nuzzrho=0}. This is achieved by now imposing
\begin{gather}
\label{newdef1}
f_1:= -\frac{1}{2}\varepsilon^2 \lambda \dot\lambda f_0g_0^{(1)}
(\frac{1}{2}) + \varepsilon^2 \lambda ^2 {\dot f}_0 g_0(\frac{1}{2}),
\\
\label{newdeff2}
f_2 :=
-\left[  (2(N-1) f_1 +\frac{1}{2}\varepsilon^2\lambda \dot \lambda)
+\frac{1}{2} \varepsilon^2 \lambda \dot\lambda
f_0 g_0 ^{(2)} (\frac{1}{2})
+ (\varepsilon^2 \lambda \dot\lambda f_0 -\varepsilon^2 \lambda ^2 {\dot f}_0)
g_0^{(1)} (\frac{1}{2}) \right],
\\
\label{newdeff3}
\begin{array}{rcl}
f_3&:=&-\left [
(2(N-1)+\frac{1}{2}\varepsilon^2\lambda \dot\lambda )f_2
+(2\varepsilon^2\lambda \dot\lambda -8(N-1))f_1 -\varepsilon^2\lambda ^2 {\dot f}_1\right.\\
&&\qquad
\left. +\frac{1}{2} \varepsilon^2\lambda \dot\lambda f_0g_0^{(3)} (\frac{1}{2})
+(2\varepsilon^2\lambda \dot\lambda f_0 -\varepsilon^2 \lambda ^2 {\dot f}_0)
g_0^{(2)} (\frac{1}{2}) \right],
\end{array}
\end{gather}
where the $g_i$'s now satisfy
\begin{gather}
\label{condiongi=1}
g_1^{(2)}\left(\frac{1}{2}\right)=g_2^{(3)}\left(\frac{1}{2}\right)=g_3^{(4)}\left(\frac{1}{2}\right)=1,
\\
\label{condiongi=0}
g_i^{(j)}\left(\frac{1}{2}\right)=0, \, \forall (i,j)\in \{1,2,3\}\times\{0,1,2,3,4\}\setminus \{(1,2),(2,3),(3,4)\}.
\end{gather}
Then $a$ still satisfies \eqref{areguliersurSigma-cas1} for some
function $\psi^a$ of class  $ C^\infty$ on $ [-1,1]\times [\rho-\delta,\rho+\delta]$. Proceeding as we did to prove
\eqref{creguliersurSigma-cas1}, we get the existence of $\psi^b$ of class  $ C^\infty$ on $ [-1,1]\times [\rho-\delta,\rho+\delta]$
such that \eqref{breguliersurSigma-cas1} holds. Now the case of the function $c$ is simpler than before since, at least for
$\varepsilon>0$ small enough, we get from \eqref{C>0pourdeltassezpetit}  that $c>0$ in $\Sigma$ and the existence
$\psi^c$ of class  $ C^\infty$ on $ [-1,1]\times [\rho-\delta,\rho+\delta]$
such that \eqref{creguliersurSigma-cas1} holds.

The case where
\begin{gather}\label{Sigmaprime}
(t,r)\in \Sigma':=\left\{(t,r)\in (-1,1)\times \R; \, z \in (-1,1)\setminus \left(\cup_{l=-1}^{p} [\rho_l-(\delta/2), \rho_l+(\delta/2)]\right)\right\}
\end{gather}
is even simpler than the two previous ones since, by \eqref{supportgil}
\begin{gather}
g_1=g_2=g_3=0.
\end{gather}
One gets that \eqref{areguliersurSigma-cas1}, \eqref{breguliersurSigma-cas1}, and \eqref{creguliersurSigma-cas1} hold on
$\Sigma'$ where
\begin{equation}
\psi^a , \, \psi^b,  \psi^c \in
C^\infty\left([-1,1]\times \left([-1,1]\setminus
\left(\cup_{l=-1}^{p} (\rho_l-(\delta/2), \rho_l+(\delta/2))\right)\right)\right).
\end{equation}
In conclusion, from these three cases we get the existence of three functions
$\psi^a$, $\psi^b$,  and $\psi^c$ such that
\begin{gather}\label{psiiregular}
\psi^a , \, \psi^b,  \psi^c \in
C^\infty([-1,1]\times [-1,1]),
\\
a(t,r)=f_0(t)\psi^a(t,z), \forall (t,r)\in \mathbb{D},
\\
b(t,r)=\lambda(t)^{-2/3}f_0(t)\psi^b(t,z), \forall (t,r) \in \mathbb{D},
\\
c(t,r)=\lambda(t)^{-8/9}f_0(t)\psi^c(t,z), \forall (t,r)\in \mathbb{D},
\end{gather}
which, together with \eqref{deflambda} and \eqref{deff0}, imply that, if $a$, $b$, and $c$ are extended to
all of $\R\times\R$ by $0$
outside $\mathbb{D}$, then $a$, $b$, and $c$ are of class $C^\infty$ on $\R\times\R$. This concludes the proof of Theorem~\ref{thmexistenceabc}. \cqfd

\section{Proof of Theorem~\ref{thnullcont}}
\label{sec-proof-local-controllability}
\setcounter{equation}{0}
In this section, we show how to deduce
Theorem~\ref{thnullcont} from Theorem~\ref{thmexistenceabc} by means
of the return method, an algebraic solvability  and classical controllability results.

Let $x_0\in \omega$. Let $\bar r >0$ be small enough so that

\begin{equation}\label{barpetit}
\left(\left|t-\frac{T}{2}\right|\leqslant \bar r^2  \text{ and }\left|x-\bar x_0\right|\leq \bar r \right)
\Rightarrow \left( t\in (0,T) \text{ and } x \in \omega\right).
\end{equation}
Let $\bar \alpha : \R\times \R^N\rightarrow \R$, $\bar \beta : \R\times \R^N\rightarrow \R$,
$\bar \gamma : \R\times \R^N\rightarrow \R$ and $\bar u : \R\times \R^N\rightarrow \R$ be defined by, for
every $(t,x)\in \R\times \R^N$,
\begin{gather}
\label{defbaralpha}
\bar \alpha(t,x):=\bar r ^8 a\left(\frac{t-(T/2)}{\bar r^2}, \frac{1}{\bar r}\left|x-x_0\right|\right),
\\
\label{defbarbeta}
\bar \beta(t,x):=\bar r ^2 b\left(\frac{t-(T/2)}{\bar r^2}, \frac{1}{\bar r}\left|x-x_0\right|\right),
\\
\label{defbargamma}
\bar \gamma(t,x):= c\left(\frac{t-(T/2)}{\bar r^2}, \frac{1}{\bar r}\left|x-x_0\right|\right),
\\
\label{defbaru}
\bar u(t,x):= \bar \gamma_t(t,x)- \Delta \bar \gamma(t,x).
\end{gather}
 From \eqref{supportabc}, \eqref{abceven}, \eqref{edpab3}, \eqref{edpbc3}, \eqref{defbaralpha}, \eqref{defbarbeta}, \eqref{defbargamma}, and \eqref{defbaru}, the functions $\bar \alpha$, $\bar \beta$, $\bar \gamma$, and
$\bar u$ are of class $C^\infty$ and satisfy
\begin{gather}
\label{eqbaralpha}
\bar \alpha_t-\Delta \bar \alpha =\bar \beta^3  \text{ in } \R \times  \R^N,
\\
\label{eqbarbeta}
\bar \beta_t-\Delta \bar \beta =\bar \gamma^3  \text{ in }\R \times  \R^N,
\\
\label{eqbargamma}
\bar \gamma_t-\Delta \bar \gamma= \bar u\chi_\omega \text{ in }\R \times  \R^N,
\\
\label{suportbarabcu}
\text{ the supports of $\bar \alpha$, $\bar \beta$, $\bar \gamma$, and $\bar u$ are included in
$(0,T) \times \omega$.}
\end{gather}
 Let $(\alpha^0,\beta^0,\gamma^0)\tr \in L^\infty(\Omega)^3$.
 For $(\alpha,\beta,\gamma)\tr\in L^\infty((0,T)\times\Omega)^3$ and $u\in L^\infty((0,T)\times\Omega)$,
let us define $(\hat \alpha,\hat \beta,\hat \gamma)\tr\in L^\infty((0,T)\times\Omega)^3$ and $\hat u\in L^\infty((0,T)\times\Omega)$ by, for every $(t,x)\in (0,T)\times \Omega$,
\begin{gather}
\label{defhatalpha}
\hat \alpha(t,x):=\alpha(t,x)-\bar \alpha(t,x),
\\
\label{defhatbeta}
\hat \beta(t,x):=\beta(t,x)-\bar \beta(t,x),
\\
\label{defhatgamma}
\hat \gamma(t,x):=\gamma(t,x)-\bar \gamma(t,x),
\\
\label{defhatu}
\hat u(t,x):=u(t,x)-\bar u(t,x).
\end{gather}
 From \eqref{eqbaralpha}, \eqref{eqbarbeta}, \eqref{eqbargamma}, and \eqref{suportbarabcu},  $(\alpha,\beta,\gamma)\tr\in L^\infty((0,T)\times\Omega)^3$  is
 the solution of the Cauchy problem \eqref{Cauchy} if and only if  $(\hat \alpha,\hat \beta,\hat \gamma)\tr\in L^\infty((0,T)\times\Omega)^3$  is
 the solution of the Cauchy problem
\begin{equation}
\label{hatCauchy}
\left\{
\begin{array}{ll}
\hat \alpha_t-\Delta \hat \alpha =3 \bar \beta ^2 \hat \beta
+  3 \bar \beta  \hat \beta^2 + \hat \beta^3 & \text{ in }(0,T)\times  \Omega,
\\
\hat \beta_t-\Delta \hat \beta =3 \bar \gamma ^2 \hat \gamma
+  3 \bar \gamma  \hat \gamma^2 + \hat \gamma^3 & \text{ in }(0,T)\times  \Omega,
\\
\hat\gamma_t-\Delta \hat\gamma= \hat u\chi_\omega& \text{ in }(0,T)\times  \Omega,
\\
\hat \alpha=\hat \beta=\hat \gamma=0 &\text{ in }(0,T)\times \partial \Omega,
\\
\hat \alpha(0,\cdot)=\alpha^0(\cdot), \, \hat \beta(0,\cdot)=\beta^0(\cdot), \,
\hat \gamma(0,\cdot)=\gamma^0(\cdot) & \text{ in } \Omega.
\end{array}
\right.
\end{equation}
Moreover, by \eqref{suportbarabcu}, \eqref{defhatalpha}, \eqref{defhatbeta}
and \eqref{defhatgamma}, one has
\begin{equation}
\label{alphahatalpha}
\alpha(T,\cdot)=\hat \alpha(T,\cdot), \, \beta(T,\cdot)=\hat \beta(T,\cdot),
\, \gamma(T,\cdot)=\hat \gamma(T,\cdot)  \text{ in }\Omega.
\end{equation}
Let us consider the system
\begin{equation}
\label{systemhat}
\left\{
\begin{array}{ll}
\hat \alpha_t-\Delta \hat \alpha =3 \bar \beta ^2 \hat \beta
+  3 \bar \beta  \hat \beta^2 + \hat \beta^3 & \text{ in }(0,T)\times  \Omega,
\\
\hat \beta_t-\Delta \hat \beta =3 \bar \gamma^2 \hat \gamma
+  3 \bar \gamma  \hat \gamma^2 + \hat \gamma^3 & \text{ in }(0,T)\times  \Omega,
\\
\hat\gamma_t-\Delta \hat\gamma= \hat u\chi_\omega& \text{ in }(0,T)\times  \Omega,
\\
\hat \alpha=\hat \beta=\hat \gamma=0 &\text{ in }(0,T)\times \partial \Omega,
\end{array}
\right.
\end{equation}
as a control system where, at time $t\in [0,T]$,  the state is $(\hat\alpha(t,\cdot), \hat \beta(t,\cdot), \hat \gamma(t,\cdot))\tr\in L^\infty(\Omega)^3$, and the control is $\hat u(t,\cdot) \in L^\infty(\Omega)$.
Note that $
(\hat\alpha, \hat \beta, \hat \gamma)\tr=0$ and $\hat u=0$ is a trajectory (i.e. a solution) of this control system.
The linearized control system around this (null) trajectory is the linear control system
\begin{equation}
\label{systemhatlinear}
\left\{
\begin{array}{ll}
\hat \alpha_t-\Delta \hat \alpha =3 \bar \beta ^2 \hat \beta
 & \text{ in }(0,T)\times  \Omega,
\\
\hat \beta_t-\Delta \hat \beta =3 \bar \gamma ^2 \hat \gamma
& \text{ in }(0,T)\times  \Omega,
\\
\hat \gamma_t-\Delta \hat \gamma= \hat u\chi_\omega& \text{ in }(0,T)\times  \Omega,
\\
\hat \alpha=\hat \beta=\hat \gamma=0 &\text{ in }(0,T)\times \partial \Omega,
\end{array}
\right.
\end{equation}
where, at time $t\in [0,T]$,  the state is $(\hat\alpha(t,\cdot), \hat \beta(t,\cdot), \hat \gamma(t,\cdot))\tr\in L^\infty(\Omega)^3$, and the control is $\hat u(t,\cdot) \in L^\infty(\Omega)$.

By \eqref{nonzeroensemble}, \eqref{defbarbeta}, and \eqref{defbargamma}, there exists  a nonempty open subset $\omega_1$ of $\omega$, $t_1\in (0,T)$ and $t_2\in (0,T)$  such that
\begin{gather}
\label{omega1etomega}
\overline {\omega_1}\subset \omega,
\\
\label{t1t3}
0<t_1<t_2<T,
\\
\label{betanot0}
\bar \beta(t,x)\not =0, \, \forall (t,x)\in [t_1,t_2]\times \overline {\omega_1},
\\
\label{gammanot0}
\bar \gamma(t,x)\not =0, \, \forall (t,x)\in [t_1,t_2]\times \overline {\omega_1}.
\end{gather}
Let $\omega_2$ be a nonempty open subset of $\omega_1$  such that
\begin{equation}\label{omega2petit}
\overline{\omega_2}\subset \omega_1.
\end{equation}
Let us
recall that, by (the proof of) \cite[Theorem 2.4, Chapter 1]{1996-Fursikov-Imanuvilov-book}, the linear control system
\begin{equation}
\label{systemhatlinear-3controls}
\left\{
\begin{array}{ll}
\hat \alpha_t-\Delta \hat \alpha =3 \bar \beta ^2 \hat \beta + v_1\chi_{(t_1,t_2)\times \omega_2}
 & \text{ in }(0,t_2)\times  \Omega,
\\
\hat \beta_t-\Delta \hat \beta =3 \bar \gamma ^2 \hat \gamma  + v_2\chi_{(t_1,t_2)\times \omega_2}
& \text{ in }(0,t_2)\times  \Omega,
\\
\hat \gamma_t-\Delta \hat \gamma= v_{3}\chi_{(t_1,t_2)\times \omega_2} & \text{ in }(0,t_2)\times  \Omega,
\\
\hat \alpha=\hat \beta=\hat \gamma=0 &\text{ in }(0,t_2)\times \partial \Omega,
\end{array}
\right.
\end{equation}
where, at time $t\in [0,t_2]$,  the state is
$(\hat\alpha(t,\cdot), \hat \beta(t,\cdot), \hat \gamma(t,\cdot))\tr\in L^\infty(\Omega)^3$
and the control is $(v_1(t,\cdot),v_2(t,\cdot),v_3(t,\cdot))\tr \in L^\infty(\Omega)^3$
 is null controllable. We next point out that, with the terminology of \cite[page 148]{1986-Gromov-book}
 (see also \cite{2014-Coron-Lissy-IM}),  the underdetermined system

\begin{equation}
\label{underdetermined}
\left\{
\begin{array}{ll}
\tilde \alpha_t-\Delta \tilde \alpha =3 \bar \beta ^2 \tilde \beta + v_1
 & \text{ in }(t_1,t_2)\times  \omega_1,
\\
\tilde \beta_t-\Delta \tilde \beta =3 \bar \gamma ^2 \tilde \gamma  + v_2
& \text{ in }(t_1,t_2)\times  \omega_1,
\\
\tilde \gamma_t-\Delta \tilde \gamma= v_3 +\tilde u  & \text{ in }(t_1,t_2)\times  \omega_1,
\end{array}
\right.
\end{equation}
where the data is $(v_1,v_2,v_3)\tr : (t_1,t_2)\times  \omega_1 \to \R^3$ and the unknown is 
$(\tilde \alpha,\tilde \beta,\tilde \gamma, \tilde u)\tr : (t_1,t_2)\times  \omega_1 \to \R^4$ is algebraically
solvable, i.e. there are solutions of \eqref{underdetermined} such that the unknown can be expressed in terms of the
derivatives of the data. Indeed, for $(v_1,v_2,v_3)\tr \in \mathcal{D}'((t_1,t_2)\times  \omega_1)^3$, if
$(\tilde \alpha,\tilde \beta,\tilde \gamma, \tilde u)\tr \in \mathcal{D}'((t_1,t_2)\times  \omega_1)^4$ is defined by
\begin{gather}
\label{defhatalpha-algebraic}
\tilde \alpha :=0,
\\
\label{defhatbeta-algebraic}
\tilde \beta :=
\displaystyle
-\frac{v_1}{3 \bar \beta ^2},
\\
\label{defhatgamma-algebraic}
\tilde \gamma :=
\displaystyle
\frac{1}{3 \bar \gamma ^2}\left(-\left(\frac{v_1}{3 \bar \beta ^2}\right)_t
+\Delta \left(\frac{v_1}{3 \bar \beta ^2}\right)-v_2 \right),
\end{gather}
\begin{equation}
\label{defhatu-algebraic}
\begin{array}{rcl}
\tilde u &:=&
\displaystyle
-v_3+ \left(\frac{1}{3 \bar \gamma ^2}\left(-\left(\frac{v_1}{3 \bar \beta ^2}\right)_t
+\Delta \left(\frac{v_1}{3 \bar \beta ^2}\right) -v_2\right)\right)_t
-\Delta\left( \frac{1}{3 \bar \gamma ^2}\left(-\left(\frac{v_1}{3 \bar \beta ^2}\right)_t
\right.\right.\\
&&
\displaystyle
\left.
\left.
 +\Delta \left(\frac{v_1}{3 \bar \beta ^2}\right) -v_2\right)\right),
\end{array}
\end{equation}
then \eqref{underdetermined} holds. This algebraic solvability is a key ingredient for the following proposition.
\begin{prop}
\label{prop-local}
There exists $\eta>0$ such that, for every $(\alpha^0,\beta^0,\gamma^0)\tr \in L^\infty(\Omega)^3$ satisfying
\begin{equation}\label{initpetit}
|\alpha^0|_{L^\infty(\Omega)}+|\beta^0|_{L^\infty(\Omega)} +|\gamma^0|_{L^\infty(\Omega)}<\eta,
\end{equation}
there exists $\hat u \in L^\infty((0,t_2)\times\Omega)$ such that the solution
$(\hat \alpha,\hat \beta,\hat \gamma)\tr \in L^\infty((0,t_2)\times \Omega)^3$ of the Cauchy problem
\begin{equation}
\label{systemhat-prop}
\left\{
\begin{array}{ll}
\hat \alpha_t-\Delta \hat \alpha =3 \bar \beta ^2 \hat \beta
+  3 \bar \beta  \hat \beta^2 + \hat \beta^3 & \text{ in }(0,t_2)\times  \Omega,
\\
\hat \beta_t-\Delta \hat \beta =3 \bar \gamma ^2 \hat \gamma
+  3 \bar \gamma  \hat \gamma^2 + \hat \gamma^3 & \text{ in }(0,t_2)\times  \Omega,
\\
\hat \gamma_t-\Delta \hat \gamma= \hat u\chi_\omega& \text{ in }(0,t_2)\times  \Omega,
\\
\hat \alpha=\hat \beta=\hat \gamma=0 &\text{ in }(0,t_2)\times \partial \Omega,
\\
\hat \alpha(0,\cdot)=\alpha^0(\cdot), \, \hat \beta(0,\cdot)=\beta^0(\cdot), \,
\hat \gamma(0,\cdot)=\gamma^0(\cdot) & \text{ in } \Omega,
\end{array}
\right.
\end{equation}
satisfies
\begin{equation}\label{valfinal}
\hat \alpha(t_2,\cdot)=\hat \beta(t_2,\cdot)=\hat \gamma(t_2,\cdot)=0 \text{ in }\Omega.
\end{equation}
\end{prop}

The proof of  Proposition~\ref{prop-local} is given in Appendix~\ref{appendix-A}.
It is an adaptation of \cite{2014-Coron-Lissy-IM}, which deals
with Navier-Stokes equations, to our parabolic system. Besides a suitable inverse mapping theorem, it mainly consists
of the following two steps.
\begin{enumerate}
\item Prove that the control system \eqref{systemhat-prop} with two ``fictitious'' controls
added on the first two equations is null controllable by means of smooth controls. See
Proposition~\ref{cor-regular-controls}.
\item Remove the two ``fictitious'' controls by using the algebraic solvability, as in \cite{1992-Coron-MCSS} and
\cite{2014-Coron-Lissy-IM}. See (the proof of) Proposition~\ref{propositionFprimeonto}.
\end{enumerate}

With the notations of Proposition \ref{prop-local}, we extend $(\hat \alpha,\hat \beta,\hat \gamma)\tr$ and $\hat u$ to all of $(0,T)\times \Omega$ by requiring
\begin{gather}\label{estensionhatabcu}
\hat \alpha (t,x)=\hat \beta (t,x)=\hat \gamma (t,x) =\hat u(t,x)=0, \, \forall (t,x) \in (t_2,T)\times \Omega.
\end{gather}
Then, by \eqref{systemhat-prop} and \eqref{valfinal}, one has \eqref{hatCauchy} and
\begin{equation}\label{valfinalhatbon}
\hat \alpha (T,\cdot)=\hat \beta (T,\cdot)=\hat \gamma (T,\cdot) =0 \text{ in }\Omega.
\end{equation}
Let us define $(\alpha,\beta,\gamma)\tr\in L^\infty((0,T)\times\Omega)^3$ and $u\in L^\infty((0,T)\times\Omega)$
by imposing \eqref{defhatalpha}, \eqref{defhatbeta}, \eqref{defhatgamma}, and \eqref{defhatu}. Then, from
\eqref{hatCauchy}, one has \eqref{Cauchy} and, using \eqref{alphahatalpha} together with \eqref{valfinalhatbon}, one has \eqref{finalnul}. This concludes the proof of Theorem~\ref{thnullcont} if \eqref{initpetit} holds.

However, assumption \eqref{initpetit} can be removed
by using the following simple homogeneity argument: If $((\alpha,\beta,\gamma)\tr, u)
\in L^\infty((0,T)\times\Omega)^3\times L^\infty((0,T)\times\Omega)$ is a trajectory
(i.e. a solution) of the control system \eqref{eqsys}, then for every $s>0$,
$((\alpha^s,\beta^s,\gamma^s)\tr, u^s):= ((s^9\alpha,s^3\beta,s\gamma)\tr, su)$ is a trajectory
(i.e. a solution) of the control system \eqref{eqsys}. This concludes the proof of Theorem~\ref{thnullcont}. \cqfd
\appendix

\section{Proof of Proposition \ref{prop-local}}
\label{appendix-A}

Let $\hat 1_{\omega_2}: \R^3\rightarrow [0,1]$ be a function of class $C^\infty$
which is equal to $1$ on $\omega_2$ and whose support is included in $\omega_1$, and
let $\zeta :\R \rightarrow [0,1]$ be such that $\zeta$ is equal to
$0$ on $(-\infty, (2t_1+t_2)/3]$ and is equal to $1$ on $((t_1+2t_2)/3, +\infty)$.
Let $\vartheta: \R\times\R^3\to \R $ be defined by
\begin{equation}\label{defvartheta}
\vartheta(t,x):=\zeta(t)\hat 1_{\omega_2}(x), \, \forall (t,x)\in \R\times\R^3.
\end{equation}
From now on, we set, $Q:=(t_1,t_2)\times \Omega$ and, for $\eta\in (0,1)$ and $K>0$,
\begin{gather}
\rho_\eta(t):=e^{\frac{-K}{\eta(t_2-t)}}, \, \rho_1(t):=e^{\frac{-K}{(t_2-t)}}, \, \forall t\in [t_1,t_2).
\end{gather}
We have the following Carleman estimates proven in \cite[Chapter 1]{1996-Fursikov-Imanuvilov-book}.
\begin{lem} \label{Carl1} Let $\eta\in (0,1)$.
There exist $K:=K(\eta)>0$ and $C:=C(K)>0$ such that, for every $g=(g_1,g_2,g_3)\tr \in L^2((t_1,t_2)\times\Omega)^3$ and for every  solution $z=(\hat \alpha, \hat \beta, \hat \gamma)\tr \in L^2((t_1,t_2),H^2(\Omega)^3)\cap H^1((t_1,t_2),L^2(\Omega)^3)$ of the parabolic system, which is the adjoint of
\eqref{systemhatlinear},
\begin{equation}
\label{Sbis-adj}
\left\{
\begin{array}{ll}
-\hat \alpha_t-\Delta \hat \alpha =g_1
 & \text{ in }(t_1,t_2)\times  \Omega,
\\
-\hat \beta_t-\Delta \hat \beta -3 \bar \beta ^2 \hat \alpha=g_2
& \text{ in }(t_1,t_2)\times  \Omega,
\\
-\hat \gamma_t-\Delta \hat \gamma- 3 \bar \gamma ^2 \hat \beta=g_3 & \text{ in }(t_1,t_2)\times  \Omega,
\\
\hat \alpha=\hat \beta=\hat \gamma=0 &\text{ in } (t_1,t_2)\times \partial \Omega,
\end{array}
\right.
\end{equation}
 one has
\begin{equation}
|\sqrt{\rho_\eta}z|^2_{L^2(Q)^3}+|z(t_1,\cdot )|_{L^2(\Omega)^3}^2\\\leqslant C\left(\int_{(t_1,t_2)\times\Omega}\vartheta  \rho_1 |z|^2+\int_{(t_1,t_2)\times\Omega}\rho_1|g|^2\right ) .\label{Carleman-good}
\end{equation}
\end{lem}

Let us now derive from Lemma~\ref{Carl1} a proposition on the null-controllability with controls which are smooth functions for the control system \eqref{systemhatlinear} with a right hand side term.

\begin{prop}
\label{cor-regular-controls}
Let $\eta \in (0,1)$ be such that
\begin{equation}\label{eta-proche-de-1-1}
\eta >\frac{2}{3}
\end{equation}
and let $K$ be as in Lemma~\ref{Carl1}. Let
 $k\in \mathbb{N}$ and let $p\in [2,+\infty)$. Then, for every $f=(f_1,f_2,f_3)\tr \in L^p(Q)^3$
   such that ${\rho_\eta}^{-1/2}f\in L^p(Q)^3$ and for every
  $(\alpha^0,\beta^0,\gamma^0)\tr \in W^{1,p}_0(\Omega)^3\cap W^{2,p}(\Omega)^3$,
  there exists $u=(u_1,u_2,u_3) \in L^2(Q)^3$ satisfying
\begin{gather}
\label{bonne-reg-control}
e^{\frac{K\eta^2}{2(t_2-t)}}\vartheta u  \in L^2((t_1,t_2),H^{2k}(\Omega)^3)\cap H^{k}((t_1,t_2),L^2(\Omega)^3),
\end{gather}
such that the solution $\hat y:=(\hat \alpha ,\hat \beta,  \hat \gamma)\tr$ of
\begin{equation}
\label{systemhat-prop-ui}
\left\{
\begin{array}{ll}
\hat \alpha_t-\Delta \hat \alpha =3 \bar \beta ^2 \hat \beta
 + f_1+ \vartheta  u_1& \text{ in }(t_1,t_2)\times  \Omega,
\\
\hat \beta_t-\Delta \hat \beta =3 \bar \gamma ^2 \hat \gamma
+ f_2+ \vartheta u_2& \text{ in }(t_1,t_2)\times  \Omega,
\\
\hat \gamma_t-\Delta \hat \gamma= f_3+ \vartheta u_3 & \text{ in }(t_1,t_2)\times  \Omega,
\\
\hat \alpha=\hat \beta=\hat \gamma=0 &\text{ in }(t_1,t_2)\times \partial \Omega,
\\
\hat \alpha(t_1,\cdot)=\alpha^0(\cdot), \, \hat \beta(t_1,\cdot)=\beta^0(\cdot), \,
\hat \gamma(t_1,\cdot)=\gamma^0(\cdot) & \text{ in } \Omega,
\end{array}
\right.
\end{equation}
satisfies
\begin{gather}
\label{estimateyLp}
e^{\frac{K}{2(t_2-t)}}\hat y \in L^p((t_1,t_2),W^{2,p}(\Omega)^3)\cap W^{1,p}((t_1,t_2),L^p(\Omega)^3).
\end{gather}
\end{prop}
\textbf{Proof of Proposition \ref{cor-regular-controls}.} We adapt the proof
of \cite[Proposition 4]{2014-Coron-Lissy-IM}
to our situation. Modifying if necessary $f$, we may assume without loss of generality that
\begin{equation}\label{alpha0=0}
  (\alpha^0,\beta^0,\gamma^0)\tr=0.
\end{equation}
Let us define a linear operator $S: \mathcal{D}'(Q)^3\rightarrow \mathcal{D}'(Q)^3$ by
\begin{gather}
\label{defSdroit}
Sz:=
\begin{pmatrix}
-\alpha_t-\Delta \alpha
\\
-\beta_t-\Delta \beta -3\bar \beta ^2 \alpha
\\
-\gamma_t-\Delta \gamma -3\bar \gamma^2 \beta
\end{pmatrix},\,
\forall z=
\begin{pmatrix}
\alpha
\\
\beta
\\
\gamma
\end{pmatrix}
\in \mathcal{D}'(Q)^3.
\end{gather}
We define a closed linear unbounded operator $\mathcal{S}: \mathcal D(\mathcal{S}) \subset L^2(Q)^3\rightarrow L^2(Q)^3$ by
\begin{multline}
\label{domainecalN}
\mathcal D(\mathcal{S}):=\{z=(\alpha,\beta,\gamma)\tr \in L^2((t_1,t_2),H^1_0\cap H^2(\Omega)^3)
\\
\cap H^1((t_1,t_2),L^2(\Omega)^3);\, z(t_2,\cdot)=0\},
\end{multline}
\begin{gather}
\mathcal{S}z=Sz.
\end{gather}
Let
\begin{equation}\label{defX0}
X_0:=L^2(Q).
\end{equation}
 For $m\in \mathbb{N}\setminus\{0\}$, we set
\begin{gather}
\label{defXm}
X_m:=\mathcal D({\mathcal S}^m),
\end{gather}
Let us point out that
\begin{equation}\label{defscalarm}
<z_1,z_2>_{X_m}:=<\mathcal S ^mz_1,\mathcal S^mz_2>_{L^2(Q)^3}
\end{equation}
is a scalar product on $X_m$. From now on  $X_m$ is equipped with this scalar product. Then $X_m$ is an Hilbert space. For $m\in \mathbb{Z}\cap (-\infty,0)$, let
\begin{gather}
\label{defX-m}
X_m:=X_{-m}',
  \end{gather}
where $X_{-m}'$ denotes the dual space of $X_{-m}$. We choose the
 pivot space $L^2(Q)^3=X_0$. In particular \eqref{defX-m} is an equality for $m=0$.
 For every $(k,l)\in \mathbb Z^2$ such that $k\leqslant l$, one has
\begin{equation}\label{XlsubsetXk}
X_l\subset X_k.
\end{equation}
Note that, since $\Omega$ is only of class $C^2$, in general, for $m\in \mathbb{N}\setminus\{0,1\}$,
\begin{equation}\label{not-subset}
X_m \not \subset L^2((t_1,t_2),H^{2m}( \Omega)^3)\cap H^{m}((t_1,t_2),L^2(\Omega)^3).
\end{equation}
However, even with  $\Omega$  only of class $C^2$, by classical results
on the interior regularity of parabolic systems, for every $m\in \mathbb{N}$, for
every open subset $\Omega_0$ such that $\overline{\Omega_0}\subset \Omega$, and for every $z\in  X_m$,
\begin{equation}\label{regularity-Xm}
z_{|(t_1,t_2)\times\Omega_0}\in L^2((t_1,t_2),H^{2m}( \Omega_0)^3)\cap H^{m}((t_1,t_2),L^2(\Omega_0)^3).
\end{equation}
(Note that this property is not known to hold for the linearized Navier-Stokes equations considered in
\cite{2014-Coron-Lissy-IM} for $\Omega$ only of class $C^2$ ; this is why $\Omega$ is assumed to be of class $C^\infty$
in \cite{2014-Coron-Lissy-IM}.)

 For $m\in\mathbb N$, one can define $\mathcal{S}^*$ as an operator from $X_{-m}$ into $X_{-m-1}$ by setting, for every $z_1\in X_{-m-1}$ and for every $z_2\in X_{m+1}$,
\begin{gather}
\label{defdualL}<\mathcal{S}^*z_1,z_2>_{X_{-m-1},X_{m+1}}:=<z_1,\mathcal S z_2>_{X_{-m},X_m}.\end{gather}
(One easily checks that this definition is consistent: it gives the same image if $z_1$ is also in $X_{-m'}$ for some $m'\in \mathbb{N})$.
This implies in particular that, for every $z_1\in L^2(Q)^3$ and for every $z_2\in X_m$, one has, for every $0\leqslant j\leqslant l$,
\begin{gather}\label{dualLkl}<({\mathcal S}^*)^lz_1,z_2>_{X_{-l},X_{l}}=<({\mathcal S}^*)^{l-j} z_1,(\mathcal S)^jz_2>_{X_{j-l},X_{l-j}}.
\end{gather}
Let $\mathcal{H}_0$ be the set of
 $z\in H^1((t_1,t_2), L^2(\Omega)^3)\cap L^2((t_1,t_2),H^2(\Omega)^3)$  such that
\begin{gather}
\label{condition1}
\sqrt{\rho_1} S z\in X_{k},
\\
\label{condition2}
\sqrt{\vartheta \rho_1} z \in L^2(Q)^3.
\end{gather}
Let $q$ be the following bilinear form defined on $\mathcal{H}_0$:
\begin{gather}
\label{defquadratic}
q(z,w):=<\sqrt{\rho_1}S z,\sqrt{\rho_1}S w>_{X_{k}}+\int_{Q}\vartheta \rho_1 z\cdot w.
\end{gather}
 (This is the analogue of the bilinear form denoted by $a$ in \cite{2014-Coron-Lissy-IM}.) From  \eqref{Carleman-good}, we deduce that $q$ is a scalar product on $\mathcal{H}_0$. Let $\mathcal{H}$ be the completion of  $\mathcal{H}_0$ for this scalar product. Note that, still from \eqref{Carleman-good} and also from the definition of
 $\mathcal{H}$,  $\mathcal{H}$ is a subspace of $L^2_{loc}((t_1,t_2),H^1_0(\Omega)^3)$ and, for every $z\in \mathcal{H}$,
 one has \eqref{condition1},  \eqref{condition2}, and
\begin{gather}
\label{estzH}
|{\rho_\eta}^{1/2}z|_{L^2(Q)^3}\leqslant C\sqrt{q(z,z)}, \mbox{ }\forall z\in \mathcal{H}.
\end{gather}
As in \cite{2014-Coron-Lissy-IM}, using the Riesz representation theorem together with \eqref{estzH},
 one gets that there exists a unique
\begin{gather}
\label{hatzinH}
\hat z\in \mathcal{H}
\end{gather}
 verifying, for every
$w\in \mathcal{H}$,
\begin{gather}\label{Lax-M}
<\mathcal S^{k}(\sqrt{\rho_1}S\hat z),\mathcal S^{k}(\sqrt{\rho_1} S w)>_{L^2(Q)^3}-\int_{Q} u\cdot w=\int_{Q} f\cdot w,
\end{gather}
with
\begin{gather}
\label{defu=rho1hatz}
u:= -\rho_1 \hat z.
\end{gather}
We then set
\begin{gather}
\label{def-tildey}
\tilde y:= (\mathcal S^*)^{k}\mathcal {S}^{k}(\sqrt{\rho_1} S \hat z)\in X_{-k}.
\end{gather}
We want to gain regularity on $\tilde y$ by accepting to have a weaker exponential decay rate for $\tilde y$ when $t$ is close to $t_2$ (in the spirit of \cite[Theorem 2.4, Chapter 1]{1996-Fursikov-Imanuvilov-book} and \cite{2002-Barbu-NA}).
Let $\psi \in C^\infty([t_1,t_2])$ and $y\in X_{-1}$. One can define $\psi y\in X_{-1}$ in the following way. Since $\mathcal S^*:X_{0}\rightarrow X_{-1}$ is onto, there exists  $h\in X_0$ such that
$\mathcal S^*h=y$. We define $\psi y $ by
\begin{gather}
\psi y= \psi \mathcal S^*h:=-\psi 'h+\mathcal S^*( \psi h).
\end{gather}
This definition is compatible with the usual definition of $\psi y$ if $y\in X_0$. We can then define by induction on $m$
$\psi y\in X_{-m}$ for $\psi \in C^\infty([t_1,t_2])$ and $y\in X_{-m}$ in the same way.
Using \eqref{def-tildey}, this allows us to define
\begin{gather}
\label{defhaty}
\hat y:= \sqrt{\rho_1} \tilde y\in X_{-k}.
\end{gather}
From \eqref{Lax-M}, \eqref{defu=rho1hatz}, \eqref{def-tildey}, and \eqref{defhaty}, one gets
\begin{gather}
\label{eqhaty}
 \mathcal{S}^* \hat y= f +\vartheta u \text{ in }X_{-k-1}.
\end{gather}
Let
\begin{gather}\label{prpertytileK1-deftilderho1}
\tilde K\in(0,K) \text{ and } \tilde \rho_1:=e^{-\tilde{K}/(t_2-t)}.
\end{gather}
Using \eqref{defu=rho1hatz}, \eqref{def-tildey}, and \eqref{eqhaty}, one has
\begin{gather}
\label{eqN*y}
\mathcal{S}^* \left(\left(\sqrt {\rho_1}/\sqrt {\tilde \rho_1}\right)\tilde y\right)
=\left(1/\sqrt{\tilde {\rho_1}}\right)'\sqrt{\rho_1} \tilde y +\left(1/\sqrt {\tilde \rho_1}\right) (f+u)
\text{ in }  X_{-k}.
\end{gather}
We want to deduce from \eqref{eqN*y} some information on the regularity of $\tilde y$. This can be achieved thanks to the following lemma, the proof of which is similar to the proof of \cite[Lemma 4]{2014-Coron-Lissy-IM}.
\begin{lem}\label{reg-dual} Let $m\in \mathbb{N}$.
If $y\in X_{-m}$ and $\mathcal S^*y\in X_{-m}$, then $y\in X_{-m+1}$.
\end{lem}
 From \eqref{defhaty}, \eqref{eqN*y}, and Lemma~\ref{reg-dual}, one gets that
$$\left(\sqrt{\rho_1}/\sqrt {\tilde \rho_1}\right)\tilde y\in X_{-k+1},\,\forall \tilde K\in (0,K) .$$
Using an easy induction argument together with Lemma~\ref{reg-dual} (and the fact that one can choose $\tilde K<K$ arbitrarily close to $K$), we deduce that, for every $\tilde K\in (0,K)$,
$\left(\sqrt{\rho_1}/\sqrt {\tilde \rho_1}\right)\tilde y\in X_{0}$.

Let us now focus on $u$. Let us define
\begin{equation}\label{defv}
v:=\rho_1 \hat z.
\end{equation}
Using \eqref{estzH}, one gets that
\begin{gather}
\label{vx0}{\rho_1}^{-1}{\rho_\eta}^{1/2}v\in L^2(Q)^3.
\end{gather}
 Using \eqref{hatzinH} together with regularity results for $S$ applied on ${\tilde \rho_1}^{-1}{\rho_\eta}^{1/2}v\in L^2(Q)^3 $ and, as above for the proof of \eqref{vx0}, a bootstrap argument (together with the fact that one can choose $\tilde K\in (0,K)$ arbitrarily close to $K$), one obtains that
\begin{gather}
\label{reg-cont-ok}{\tilde\rho_1}^{-1}{\rho_\eta}^{1/2}v\in X_{k}, \, \forall \tilde K\in (0,K).
\end{gather}
Let us point out that \eqref{eta-proche-de-1-1} implies that
\begin{equation}\label{eta2}
  \eta^2-2+\frac{1}{\eta}<0.
\end{equation}
 From \eqref{eta-proche-de-1-1}, \eqref{regularity-Xm}, \eqref{defu=rho1hatz}, \eqref{defv},
\eqref{reg-cont-ok}, and \eqref{eta2}, one gets \eqref{bonne-reg-control}.

Let us now deal with $\hat y$. Without loss of generality, we may assume that
\begin{equation}\label{kassezgrand-1}
4k> 2+N,
\end{equation}
so that
\begin{equation}\label{regularity-L-infty}
L^2((t_1,t_2),H^{2k}(\Omega)^3)\cap H^{k}((t_1,t_2),L^2(\Omega)^3)\subset L^\infty(Q).
\end{equation}
From \eqref{eqhaty}, \eqref{regularity-L-infty}, and \eqref{reg-cont-ok}, we deduce (by looking at the parabolic system verified by $(1/\sqrt {\tilde \rho_1}) \hat y$ and using usual regularity results for linear parabolic systems) that
\begin{gather}
\label{reg-sol-ok}
\left(1/\sqrt {\tilde \rho_1}\right) \hat y \in L^p((t_1,t_2),W^{2,p}(\Omega)^3)\cap W^{1,p}((t_1,t_2),L^p(\Omega)^3), \, \forall \tilde K\in (0,K),
\end{gather}
which, together with \eqref{eqhaty}, concludes the proof of Proposition~\ref{cor-regular-controls}.
\cqfd

To end the proof of Proposition~\ref{prop-local}, we are going to apply the following inverse mapping theorem (see \cite[Chapter 2, Section 2.3]{1987-Alekseev-et-al-book}).
\begin{prop}
\label{prop-inverse-mapping}
Let $E$ and $F$ be two Banach spaces. Let $\mathcal F:E\rightarrow F$ be of class $C^1$ in a neighborhood of $0$. Let us assume that the operator $\mathcal F'(0)\in\mathcal L(E,F)$ is onto. Then there exist $\eta>0$ and $C>0$ such that for every $g\in F$ verifying  $|g-\mathcal F(0)|<\eta$, there exists $e\in E$ such that
\begin{enumerate}
\item $\mathcal F(e)=g$,
\item $|e|_E\leqslant C|g-\mathcal F(0)|_F$.
\end{enumerate}
\end{prop}

We now use the same technique as in
\cite[Theorem 4.2]{1996-Fursikov-Imanuvilov-book}.
 For $y:=(\alpha,\beta,\gamma)\tr \in \mathcal{D}'(Q)^3$
and for $v\in \mathcal{D}'(Q)$, one defines $\mathcal{L}(y,v)\in \mathcal{D}'(Q)^3$ by
\begin{equation}\label{defL}
\mathcal{L}(y,v):=
\begin{pmatrix}
\alpha_t-\Delta \alpha -3 \bar \beta ^2 \beta
\\
\beta_t-\Delta \beta -3 \bar \gamma ^2 \gamma
\\
\gamma_t-\Delta \gamma- v
\end{pmatrix}
.
\end{equation}

Let $\eta\in (0,1)$ and let $K=K(\eta)>0$ be as in Lemma~\ref{Carl1}.
We apply Proposition~\ref{prop-inverse-mapping} with $E$ and $F$ defined in the following way.
Let $E$ be the space of the functions
$$(y,v)\in L^p(Q)^3\times L^\infty(Q)$$ such that
\begin{enumerate}
\item $e^{\frac{\eta^3 K}{2(t_2-t)}}y \in L^p((t_1,t_2),W^{2,p}(\Omega)^3)\cap W^{1,p}((t_1,t_2),L^p(\Omega)^3)$,
\item $e^{\frac{\eta^3 K}{2 (t_2-t)}}v\in L^{\infty}(Q)^3$ and the support of $v$ is included in $(t_1,t_2)\times \omega $,
\item $e^{\frac{K}{2\eta (t_2-t)}}\mathcal{L}(y,v)\in L^p(Q)^3$,
\item $y(t_1,\cdot )\in W^{1,p}_0(\Omega)^3\cap W^{2,p}(\Omega)^3$,
\end{enumerate}
 equipped with the following norm which makes it a Banach space:
\begin{multline}
\label{def-norm-E}
|(y,v)|_E:= |e^{\frac{\eta^3 K}{2(t_2-t)}}y|_{L^p((t_1,t_2),W^{2,p}(\Omega)^3)\cap W^{1,p}((t_1,t_2),L^p(\Omega)^3)}\\
+ |e^{\frac{\eta ^3 K}{2(t_2-t)}}v|_{L^\infty(Q)}+|e^{\frac{K}{2\eta (t_2-t)}}\mathcal{L}(y,v)|_{L^p(Q)^3}
+|y(t_1,\cdot )|_{W^{2,p}(\Omega)^3}.
\end{multline}
Let $F$ be the space of the functions $(h,y^0)\in L^p(Q)^3
\times \left(W^{1,p}_0(\Omega)^3\cap W^{2,p}(\Omega)^3\right)$ such that
\begin{equation}\label{defF}
e^{\frac{K}{2\eta (t_2-t)}}h\in L^p(Q)^3
\end{equation}
equipped with the following norm which makes it a Banach space:
\begin{equation}\label{defnormsurF}
|(h,y^0)|_F:=|e^{\frac{\eta {K}_1}{2(t_2-t)}}h|_{L^p(Q)^3}+ |y^0|_{W^{2,p}(\Omega)^3}.
\end{equation}

We define $\mathcal F: E\rightarrow F$ by
\begin{equation}\label{defapplicationF}
\mathcal F(y,v)=\left(\mathcal{L}(y,v)-\begin{pmatrix}
3 \bar \beta  \beta^2 + \beta^3
\\
 3 \bar \gamma  \gamma^2 +  \gamma^3
\\
0
\end{pmatrix},y(t_1,\cdot )\right).
\end{equation}
One easily sees that $\mathcal F$ is of class $C^1$ if
\begin{equation}\label{p-largenough-pouralgebre}
  p>\frac{N+2}{2} \text{ and } \eta >\frac{1}{2^{1/4}}.
\end{equation}
 From now on, we assume  $p>2$ and $\eta \in (0,1)$ are  chosen so that \eqref{p-largenough-pouralgebre} holds.
 Note that the second inequality of \eqref{p-largenough-pouralgebre} implies that \eqref{eta-proche-de-1-1} holds. Let us assume for the moment
 that the following proposition holds.
\begin{prop}
\label{propositionFprimeonto}
One has
\begin{equation}\label{Fprimeonto}
\mathcal{F}'(0,0)(E)=F.
\end{equation}
\end{prop}
Then the assumptions of Proposition~\ref{prop-inverse-mapping} hold. Since Proposition~\ref{prop-local}
follows from the conclusion of Proposition~\ref{prop-inverse-mapping} by taking $\hat u =0$ in $(0,t_1)\times \Omega$,
this concludes the proof of
 Proposition~\ref{prop-local}.

It only remains to prove Proposition~\ref{propositionFprimeonto}. Let $f=(f_1,f_2,f_3)\tr$ and
$y^0= (\alpha^0,\beta^0,\gamma^0)\tr $ be such that $(f,y^0)\in F$.
Let us choose $k$ large enough so that
\begin{equation}\label{kassez-grand}
N+2<4(k-2).
\end{equation}
Using  Proposition~\ref{cor-regular-controls}, we get the existence of  $u=(u_1,u_2,u_3) \in L^2(Q)^3$ satisfying \eqref{bonne-reg-control}
such that the solution $\hat y:=(\hat \alpha ,\hat \beta,  \hat \gamma)\tr$ of \eqref{systemhat-prop-ui}
satisfies \eqref{estimateyLp}.
We now use the algebraic solvability of \eqref{underdetermined}
(i.e. that \eqref{defhatalpha-algebraic}, \eqref{defhatbeta-algebraic}, \eqref{defhatgamma-algebraic},
and \eqref{defhatu-algebraic}
 imply \eqref{underdetermined}) with
\begin{equation}\label{defvvarthetau}
v:=\vartheta u.
\end{equation}
We get that, if
 $(\tilde \alpha,\tilde \beta,\tilde \gamma, \tilde u)\tr \in \mathcal{D}'((t_1,t_2)\times  \omega_1)^4$
 is defined by \eqref{defhatalpha-algebraic}, \eqref{defhatbeta-algebraic}, \eqref{defhatgamma-algebraic},
 and
\eqref{defhatu-algebraic},
then \eqref{underdetermined} holds. We extend $\tilde \alpha$, $\tilde \beta$, $\tilde \gamma$,
and $\tilde u$ to $(t_1,t_2)\times \Omega$ by
$0$ outside $(t_1,t_2 )\times (\Omega \setminus \omega_1)$ and
still denote by $\tilde \alpha$, $\tilde \beta$, $\tilde \gamma$,
and $\tilde u$ these extensions. Note that \eqref{underdetermined} still holds on
$(t_1,t_2)\times \Omega$ and that (see, in particular \eqref{defvartheta})
\begin{gather}
\label{valeurt1}
\tilde \alpha(t_1,\cdot)= \tilde \beta(t_1,\cdot)= \tilde \gamma(t_1,\cdot)= 0.
\end{gather}
Finally we define $y:=(\alpha,\beta,\gamma)\tr \in \mathcal{D}'((t_1,t_2)\times  \omega_1)^3$ and
$u\in \mathcal{D}'((t_1,t_2)\times  \omega_1)$
by
\begin{equation}\label{defalphabetagammau}
\alpha:=\hat \alpha -\tilde \alpha,  \,\beta:=\hat \beta -\tilde \beta, \,\, \gamma:=\hat \gamma -\tilde \gamma,
u:=-\tilde u.
\end{equation}
 From \eqref{defhatalpha-algebraic}, \eqref{defhatbeta-algebraic}, \eqref{defhatgamma-algebraic}
\eqref{defhatu-algebraic}, \eqref{bonne-reg-control}, \eqref{estimateyLp}, \eqref{kassez-grand},
\eqref{defvvarthetau},  and \eqref{defalphabetagammau}, we get that $(y,u)\in E$. Then,
 from \eqref{underdetermined}, \eqref{systemhat-prop-ui}, \eqref{defvvarthetau}, \eqref{valeurt1}, and
 \eqref{defalphabetagammau}, we get that $\mathcal{F}'(0,0)(y,u)=(y^0,f)$.
 This concludes the proof of Proposition~\ref{propositionFprimeonto} and therefore also
 the proof of Proposition~\ref{prop-local}.
\cqfd

\begin{rem}
\label{rem-otherapproach} 1. Instead of proceeding as in \cite{2014-Coron-Lissy-IM} in order to prove Proposition~\ref{prop-local}, one can also proceed as in \cite{2010-Coron-Guerrero-Rosier-SICON}. For that, an important step is to prove that small (in a suitable sense) perturbations of the linear control system \eqref{systemhatlinear} are controllable by means of bounded controls (see \cite[Section 3.1.2]{2010-Coron-Guerrero-Rosier-SICON}. This controllability property follows from
\cite[Theorem 4.1]{2010-Gonzalez-Teresa-PM} and one can also get it by following
\cite[Section 3.1.2]{2010-Coron-Guerrero-Rosier-SICON} or \cite{2015-Duprez-Lissy-HAL}. 2. Let us emphasize that the algebraic solvability of \eqref{underdetermined} leads to a loss of derivatives. This problem is managed in our situation thanks to hypoelliptic
properties of parabolic equations. These properties do not hold, for example, for hyperbolic equations. However, for these last equations, the loss of derivatives problem can be solved thanks to a Nash-Moser inverse mapping theorem due to
Gromov \cite[Section 2.3.2, Main Theorem]{1986-Gromov-book}.  See \cite{2015-Alabau-Coron-Olive-HAL} for the first use of this inverse mapping theorem in the context of control of partial differential equations.
\end{rem}

\def\cprime{$'$}

\end{document}